\documentclass[10 pt,a4paper,leqno]{article}
\usepackage{latexsym,amsmath,amscd,amsfonts}
\usepackage[francais,english]{babel}
\usepackage[all]{xy}

\title{La filtration canonique par les pentes 
d'un module aux $q$-diff\'erences
et le gradu\'e associ\'e.}

\author{\small Jacques Sauloy, \\
\small Laboratoire Emile Picard, UMR 5580,\\
\small Universit\'e Paul Sabatier, U.F.R. M.I.G.\\
\small 118, route de Narbonne, 31062 Toulouse Cedex 4, France\\
\small {\tt sauloy@picard.ups-tlse.fr}
}

\date{\today}

\newcommand{\Ker}{\text{Ker }}

\def \gt {\symbol{62}}

\def \lt {\symbol{60}}

\def \Pr {\textsl{Preuve. - }}

\begin{document}

\maketitle

\bigskip \hrule \bigskip

\centerline{\textbf{\emph{R\'esum\'e}}}

\emph{Nous montrons que le polygone de Newton d'une \'equation 
aux $q$-diff\'erences lin\'eaire ne d\'epend que du module aux 
$q$-diff\'erences correspondant. Nous interpr\'etons les classiques 
r\'esultats de factorisation convergente de Adams-Birkhoff-Guenther 
en termes d'existence d'une filtration canonique par les pentes. 
De plus, le gradu\'e associ\'e poss\`ede d'excellentes propri\'et\'es 
fonctorielles (d'o\`u son inter\^et pour la classification) et 
tensorielles (d'o\`u son inter\^et pour la th\'eorie de Galois).}

\bigskip \hrule \bigskip

\centerline{\textbf{\emph{Abstract}}}

\emph{We show that the Newton polygon of a linear $q$-difference 
equation depends only on the corresponding $q$-difference module. 
We interpret the classical results of convergent factorisation of 
Adams-Birkhoff-Guenther in terms of the existence of a canonical 
filtration. Moreover, the associated graded module has excellent 
functorial (resp. tensorial) properties, whence its interest for 
classification (resp. for Galois theory).}

\bigskip \hrule \bigskip

\tableofcontents

\bigskip \hrule \bigskip

%%%%%%%%%%%%%%%%%%%%%%%%%%%%%%%%%%%%%%%%%%%%%%%%%%%%%%%%%%%%%%%%%%%%%%%%%%%%%%

\setcounter{equation}{-1}

\section*{Introduction}

Soit $q$ un nombre complexe de module $|q| \gt 1$.
L'\'etude locale de l'\'equation aux $q$-diff\'erences lin\'eaire:
\begin{equation}
a_{0}(z) \; f(q^{n} z) + a_{1}(z) \; f(q^{n-1} z) + 
\cdots + a_{n}(z) \; f(z) = 0,
\end{equation}
fait intervenir un polygone de Newton (\emph{voir} 
\cite{Adams1},\cite{Adams2}). Nous montrons que celui-ci est en fait un 
objet intrins\`eque: il peut \^etre d\'efini en fonction du module aux 
$q$-diff\'erences associ\'e. Il poss\`ede en outre de bonnes propri\'et\'es
fonctorielles, ab\'eliennes et tensorielles. Dans le cas formel, il donne 
lieu \`a une d\'ecomposition du module aux $q$-diff\'erences associ\'e
en somme directe de modules purs (\`a une seule pente), comme
dans le cas classique des \'equations diff\'erentielles complexes
(\emph{voir} \cite{Praagman}, \cite{SVdP}). \\

Dans le cas convergent, apparait un ph\'enom\`ene sp\'ecifique aux 
$q$-diff\'erences: le lemme d'Adams garantit l'existence de solutions 
convergentes associ\'ees \`a la premi\`ere pente. Birkhoff et Guenther 
en ont d\'eduit dans \cite{Birkhoff3} une factorisation canonique 
convergente de tout op\'erateur aux $q$-diff\'erences. Ces r\'esultats
ont \'et\'e repris, am\'elior\'es  et utilis\'es par Marotte et Zhang 
(\emph{voir} \cite{MarotteZhang}). Nous en donnons une interpr\'etation 
en termes d'existence, pour tout module aux $q$-diff\'erences, d'une 
filtration canonique par les pentes avec des quotients purs. \\

Birkhoff fondait de grands espoirs sur la factorisation canonique pour la 
formation d'invariants transcendants (\emph{voir loc. cit.}). Nous montrons 
que le foncteur ``gradu\'e associ\'e'' poss\`ede en effet d'excellentes 
propri\'et\'es fonctorielles, ab\'eliennes et tensorielles, qui ont permis
(avec d'autres outils plus puissants) d'achever le programme de Birkhoff
(\emph{voir} \cite{RSZ}). \\

Ces m\^emes propri\'et\'es, proches de celles axiomatis\'ees par
Saavedra dans \cite{Saavedra}, permettent \'egalement de passer de la
th\'eorie de Galois des \'equations fuchsiennes, d\'evelopp\'ee par
voie analytique dans \cite{JSGAL}, \`a la th\'eorie de Galois locale 
des \'equations irr\'eguli\`eres (\emph{voir} \cite{JSIRR}).

%%%%%%%%%%%%%%%%%%%%%%%%%%%%%%%%%%%%%%%%%%%%%%%%%%%%%%%%%%%%%%%%%%%%%%%%%%%%%%

\subsection*{Organisation de cet article}

Dans la premi\`ere section, nous reprenons et mettons en forme des
r\'esultats classiques dus, pour l'essentiel, \`a Adams et \`a Birkhoff.
Ces \'enonc\'es ont \'et\'e exhum\'es apr\`es un long sommeil par 
Changgui Zhang: voir \cite{MarotteZhang} et \cite{ZhangSommable}
(en particulier le paragraphe 5), o\`u l'on trouvera \'egalement des 
variantes des formulations ci-dessous. Ces r\'esultats vont par paire: 
cas convergent-cas formel. Les principaux \'enonc\'es sont les th\'eor\`emes 
de factorisation: propositions 1.2.4 et 1.2.5 et (surtout) le th\'eor\`eme 
1.2.8. \\

Dans la deuxi\`eme section, apr\`es de brefs rappels sur le formalisme
des modules aux $q$-diff\'erences, nous montrons le caract\`ere intrins\`eque 
du polygone de Newton (th\'eor\`eme 2.2.6). L'ingr\'edient principal est
l'utilisation du th\'eor\`eme de Jordan-H\"{o}lder, selon la m\'ethode de 
Katz dans \cite{Katz2}, II.2.2. Nous d\'ecrivons ensuite le comportement du
polygone de Newton vis \`a vis des op\'erations lin\'eaires. \\

L'\'etude au 3.1 du sous-module de rang maximum de pente donn\'ee
est la deuxi\`eme \'etape cruciale. Le th\'eor\`eme 3.1.1 est une 
traduction du lemme d'Adams. On en d\'eduit facilement l'existence de 
la filtration canonique (th\'eor\`eme 3.1.6), qui est une traduction 
du th\'eor\`eme de factorisation de Birkhoff-Guenther. Les excellentes 
propri\'et\'es de la filtration et du gradu\'e associ\'e  vis \`a vis 
des op\'erations lin\'eaires sont donn\'ees en 3.2 et 3.3. \\

Nous esquissons enfin en 3.4 des applications de ces r\'esultats
\`a la classification et \`a la th\'eorie de Galois; celles-ci
feront l'objet de publications ult\'erieures. \\

Nous avons report\'e dans l'appendice l'application \`a la r\'esolution 
des th\'eor\`emes de factorisation, qui n'est pas logiquement n\'ecessaire
\`a nos r\'esultats (th\'eor\`emes A.3.4 et A.4.1). Outre la description
tr\`es d\'etaill\'ee des algorithmes, la principale diff\'erence avec les 
r\'ef\'erences mentionn\'ees ci-dessus est que nous n'utilisons que des 
solutions uniformes sur $\mathbf{C}^{*}$, ce qui est important pour d'autres
parties de la th\'eorie (\emph{voir} \cite{JSGAL} et \cite{JSIRR}).

%%%%%%%%%%%%%%%%%%%%%%%%%%%%%%%%%%%%%%%%%%%%%%%%%%%%%%%%%%%%%%%%%%%%%%%%%%%%%%

\subsection*{Remerciements}

Cet article 
\footnote{Les r\'esultats pr\'esent\'es ici ont \'et\'e annonc\'es dans 
une note parue aux C.R.A.S. en janvier 2002. Le texte comportait une
petite erreur, corrig\'ee ici (section 2.2).}
provient, pour l'essentiel, de la r\'edaction d'expos\'es 
au Groupe de Travail sur les Equations aux $q$-Diff\'erences, dont je
remercie tout particuli\`erement les animateurs, Lucia Di Vizio et
Jean-Pierre Ramis, ainsi que l'un des participants \'episodiques,
Changgui Zhang, pour de nombreuses discussions excitantes, d'utiles
conseils, et le plaisir d'une passion partag\'ee. \\

\emph{Ce travail est d\'edi\'e \`a Jean Giraud, dont le cours \`a Orsay
``Etude locale des singularit\'es'' m'a appris l'inter\^et des filtrations
en g\'eom\'etrie, et, plus g\'en\'eralement, le plaisir de l'outil
bien fait.}

%%%%%%%%%%%%%%%%%%%%%%%%%%%%%%%%%%%%%%%%%%%%%%%%%%%%%%%%%%%%%%%%%%%%%%%%%%%%%%

\subsection*{Conventions g\'en\'erales}

Le corps de base $K$ est l'un des suivants:
$$
\mathcal{M}(\mathbf{C}) \subset \mathbf{C}(\{z\}) \subset \mathbf{C}((z)). 
$$
Ce qui suit s'appliquera donc en particulier aux \'equations
rationnelles, i.e. \`a coefficients dans $\mathbf{C}(z)$.
Notons que dans le cas (que nous appellerons ``classique'')
des \'equations diff\'erentielles, on ne consid\`ere pas
habituellement le corps de base $\mathcal{M}(\mathbf{C})$.
La raison pour le traiter \`a part ici est la propri\'et\'e
des \'equations aux $q$-diff\'erences \`a coefficients
rationnels de ``propager la m\'eromorphie''. \\

Le corps $K$ est muni de la valuation discr\`ete $v_{0}$ (valuation
$z$-adique). On notera $\mathcal{O}$ l'anneau de valuation correspondant,
de corps r\'esiduel $\mathcal{O}/z \mathcal{O} = \mathbf{C}$.
Le corps $K$ est \'egalement muni d'un automorphisme:
$$
\sigma_{q} : f(z) \mapsto f(qz).
$$
On consid\`erera de plus une extension $(L,\sigma_{q})$
du ``corps aux $q$-diff\'erences'' $(K,\sigma_{q})$,
o\`u l'on cherchera les solutions d'\'equations.
On supposera en particlulier que l'on peut y r\'esoudre 
les \'equations $\sigma_{q} f = z f$, $\sigma_{q} f = c f$, 
$c \in \mathbf{C}^{*}$ et $\sigma_{q} f = f + 1$.
Dans la th\'eorie des \'equations diff\'erentielles, 
ceci introduit automatiquement
des fonctions multivalu\'ees, d'o\`u la n\'ecessit\'e
d'agrandir le corps de base. Nous r\'eussirons \`a tout
faire avec des fonctions uniformes, mais $\sigma_{q}$
``propage les p\^oles'', ceux des solutions forment
des demi-spirales logarithmiques (engendr\'ees par
$q^{\pm \mathbf{N}}$) et il nous faudra tout de m\^eme
agrandir le corps de base.

\begin{itemize}

\item{Si $K = \mathcal{M}(\mathbf{C})$, on prendra en g\'en\'eral 
$L = \mathcal{M}(\mathbf{C}^{*})$.
Les \'equations \'el\'ementaires seront r\'esolues \`a l'aide 
de la fonction Theta de Jacobi, et des fonctions qui en d\'erivent.
Modifiant l\'eg\`erement les notations de \cite{JSAIF}, nous poserons 
$\Theta_{q}(z) = \theta_{q}(-z/q)$, o\`u
$\theta_{q}(z) = \underset{n \in \mathbf{Z}}{\sum}(-1)^{n}q^{-n(n-1)/2}z^{n}$.
Nous poserons \'egalement $l_{q}(z) = z\theta_{q}'(z)/\theta_{q}(z)$
et, pour tout complexe non nul $c$,
$e_{q,c}(z) = \theta_{q}(z)/\theta_{q}(c^{-1}z)$,
de sorte que $\sigma_{q}(\Theta_{q}) = z \Theta_{q}$,
$\sigma_{q}(l_{q}) = l_{q} + 1$ et $\sigma_{q}(e_{q,c}) = c e_{q,c}$.}

\item{Si $K= \mathbf{C}(\{z\})$), on prendra 
$L = \mathcal{M}(\mathbf{C}^{*},0)$, le corps
des germes en $0$ de fonctions m\'eromorphes sur $\mathbf{C}^{*}$.
Les \'equations \'el\'ementaires seront r\'esolues
\`a l'aide des germes en $0$ des fonctions pr\'ec\'edentes.}

\item{Si $K = \mathbf{C}((z))$ (``cas formel''), 
$L$ sera obtenu par adjonction de symboles permettant 
de r\'esoudre les \'equations $\sigma_{q} f = z f$,
$\sigma_{q} c = z f$ (pour tout $c \in \mathbf{C}^{*}$)
et $\sigma_{q} f = f + 1$, astreints \`a des relations
alg\'ebriques comme par exemple dans \cite{SVdP}. 
Dans ce cas, $L$ ne sera pas n\'ecessairement un corps.}

\end{itemize}
Sous ces conditions, on constate que l'on peut en fait 
r\'esoudre dans $L$ toute \'equation d'ordre $1$, avec ou sans 
second membre (respectivement: appendice ou 1.1.7). \\

Les constantes de notre th\'eorie sont les \'el\'ements
invariants par $\sigma_{q}$. Il est facile de v\'erifier
que le corps des constantes $C_{K}$ de $K$ est $\mathbf{C}$ dans
tous les cas. Celui de $L$ est encore $C_{L} = \mathbf{C}$ dans
le cas formel (\emph{voir} \cite{SVdP}). Dans le cas
dit ``convergent'' (i.e. dans les deux premiers cas),
on peut v\'erifier qu'il s'identifie au corps
$\mathcal{M}(\mathbf{E}_{q})$ des fonctions elliptiques
relatif \`a la courbe elliptique
$\mathbf{E}_{q} = \mathbf{C}^{*}/q^{\mathbf{Z}}$
(\emph{voir} \cite{JSAIF}). \\

Nous aurons besoin des extensions ramifi\'ees $K_{l}$ de $K$
pour $l \in \mathbf{N}^{*}$: elles sont d\'efinies de fa\c{c}on
naturelle \`a l'aide de variables $z_{l}$ telles que
$z_{l}^{l} = z$; on les suppose de plus compatibles, c'est \`a dire
que $z_{lm}^{l} = z_{m}$. On introduit de m\^eme une famille compatible
de racines de $q$: $q_{l}^{l} = q$ et $q_{lm}^{l} = q_{m}$. Il suffit
pour cela de fixer $\tau \in \mathbf{C}$ tel que 
$q = e^{- 2 \imath \pi \tau}$, puis de prendre
$q_{l} = e^{- 2 \imath \pi \tau/l}$. \\

On notera enfin $\mathcal{D}_{q} = K\left<\sigma,\sigma^{-1}\right>$ 
l'alg\`ebre de \"{O}re des polyn\^omes de Laurent non commutatifs, 
caract\'eris\'ee par les relations:
$$
\forall x \in K, \forall k \in \mathbf{Z} \;,\;
\sigma^{k} x = \sigma_{q}^{k}(x) \sigma^{k}.
$$
Un tel polyn\^ome $P \in \mathcal{D}_{q}$ mod\'elise donc l'op\'erateur
aux $q$-diff\'erences $P(\sigma_{q})$, d'o\`u une op\'eration 
de $\mathcal{D}_{q}$ sur $L$. On utilisera principalement des polyn\^omes 
\emph{entiers}, c'est \`a dire dont tous les mon\^omes sont \`a degr\'es 
positifs. Avec ces conventions, l'\'equation (0) s'\'ecrit:
\begin{equation}
P.f \underset{def}{=} P(\sigma_{q})(f) = 
a_{0} \; \sigma_{q}^{n} f + \cdots + a_{n} \; f = 0 \;, \quad
a_{0},\ldots,a_{n} \in K \;,\; a_{0} a_{n} \not= 0.
\end{equation}
Dans cette \'equation, l'op\'erateur aux $q$-diff\'erences $P$
est le polyn\^ome (entier) $a_{0} \; \sigma^{n} + \cdots + a_{n}$ 
de $\mathcal{D}_{q}$. La fonction inconnue $f$ est recherch\'ee
dans $L$. \\

On v\'erifie facilement que l'anneau $\mathcal{D}_{q}$ 
est euclidien (\`a gauche et \`a droite). Soit
$P = \underset{\alpha \leq i \leq \beta}{\sum} a_{i} \sigma^{i}$
un \'el\'ement de $\mathcal{D}_{q}$. On appellera 
\emph{degr\'e absolu} de $P$ l'entier naturel 
$\deg(P) = \beta - \alpha$
si $a_{\alpha} a_{\beta} \not= 0$ (et $- \infty$ si $P = 0$).
Il est imm\'ediat que $\deg(PQ) = \deg(P) + \deg(Q)$. 
On appellera \emph{valuation $z$-adique} de $P$ l'entier 
$v_{0}(P) = \min(v_{0}(a_{\alpha}),\ldots,v_{0}(a_{\beta}))$
(donc $+ \infty$ si $P = 0$). Une variante du lemme de Gauss
permet de montrer que $v_{0}(PQ) = v_{0}(P) + v_{0}(Q)$.

%%%%%%%%%%%%%%%%%%%%%%%%%%%%%%%%%%%%%%%%%%%%%%%%%%%%%%%%%%%%%%%%%%%%%%%%%%%

% 1

\section{Polygone de Newton, factorisation, solutions}

% 1.1

\subsection{Le polygone de Newton d'une \'equation aux $q$-diff\'erences}

\subsubsection*{Polygone de Newton}

Soit 
$P = \sum a_{i} \; \sigma^{i} \in \mathcal{D}_{q}$ 
un op\'erateur aux $q$-diff\'erences non nul. 
On d\'efinit son \emph{polygone de Newton} $N(P)$ comme 
l'enveloppe convexe dans $\mathbf{R}^{2}$ de l'ensemble:
$$
\{(i,j) \in \mathbf{Z}^{2} \;/\; j \geq v_{0}(a_{i}) \}.
$$
C'est aussi, par d\'efinition, le polygone de Newton de l'\'equation 
aux $q$-diff\'erences
$P.f = \sum a_{i} \; \sigma_{q}^{i}f = 0$.
On peut d'ailleurs se restreindre aux $a_{i}$ non nuls. \\

Cette d\'efinition d\'epend \'evidemment du choix de la valuation $v_{0}$.
Dans le cas d'\'equations \`a coefficients dans $K_{l}$ (obtenues par
ramification, par exemple en 1.1.4, etc ...) c'est la valuation
$z_{l}$-adique qui sera employ\'ee. \\

\textsl{1.1.1 Terminologie. -}
La fronti\`ere de $N(P)$ est form\'ee de deux demi-droites verticales
et de $k \geq 1$ vecteurs de coordonn\'ees 
$(r_{1},d_{1}),\ldots,(r_{k},d_{k}) \in \mathbf{N}^{*} \times \mathbf{Z}$,
et de pentes 
$\mu_{1} = \frac{d_{1}}{r_{1}},\ldots,\mu_{k} = \frac{d_{k}}{r_{k}}
\in \mathbf{Q}$. On suppose celles-ci rang\'ees par ordre d\'ecroissant:
$\mu_{1} > \cdots > \mu_{k}$. Les lettres $r$, $d$ sont choisies
par analogie avec des notions de rang et de degr\'e (de fibr\'es
vectoriels, par exemple). La \emph{premi\`ere pente} est $\mu_{k}$
(donc, la plus petite). \\

On notera $S(P) = \{\mu_{1},\ldots,\mu_{k}\}$ l'ensemble des pentes de $P$.
La \emph{fonction de Newton} de $P$ est la fonction
$r_{P}: \mathbf{Q} \rightarrow \mathbf{N}$ de support $S(P)$
et telle que $\mu_{i} \mapsto r_{i}$ pour $i = 1,\ldots,n$. On a donc:
$$
r_{P} = \sum_{i = 1}^{k} r_{i} \delta_{\mu_{i}},
$$
o\`u $\delta_{\mu}$ d\'esigne la fonction de Kronecker (indicatrice
de $\{\mu\}$). \\

\textsl{1.1.2 Remarque. -} On prendra garde que les d\'efinitions
ci-dessus sont adapt\'ees \`a notre convention $|q| > 1$. Pour la  m\^eme
raison, une partie des r\'esultats (voir en particulier 1.2.9 et 3.1.3)
d\'epend de l'ordre des pentes. \\

\textsl{1.1.3 Lemme. -}
\emph{La correspondance entre fonctions de Newton et polygones de Newton
est une bijection additive.} \\

\Pr
Il est facile de construire la fonction de Newton $r$ \`a partir du
polygone de Newton $N$ et r\'eciproquement. L'addition \'etant
associative et commutative des deux c\^ot\'es, il suffit, pour prouver
l'additivit\'e, de la v\'erifier dans le cas d'une somme 
$\delta_{\mu} + r$, o\`u $\mu$ minore le support de $r$. Mais, dans
ce cas, c'est un exercice facile de g\'eom\'etrie affine.
\hfill $\Box$

\subsubsection*{Manipulations \'el\'ementaires sur les pentes}

Si l'on multiplie (\`a gauche ou \`a droite) l'op\'erateur $P$
par $a \sigma^{k}$, o\`u $a \in K^{*} \;,\; k \in \mathbf{Z}$,
le polygone de Newton $N(P)$ subit une translation de vecteur 
$(k,v_{0}(a))$. En particulier, en ramenant $P$ \`a la forme (1) avec
$a_{0} = 1$, on cale l'origine du premier vecteur en $(0,0)$. 
Dor\'enavant, nous supposerons que 
$P = \sigma^{n} + a_{1} \; \sigma^{n-1} + \cdots + a_{n}$ est entier 
unitaire (par commodit\'e, nous conservons la notation $a_{0} = 1$). \\

\textsl{1.1.4 Ramification $z = z_{l}^{l} \;,\; q = q_{l}^{l}$. -}
Les pentes sont multipli\'ees par $l$. En particulier, en
prenant pour $l$ un multiple commun des $r_{i}$ (par exemple
$n !$), on se ram\`ene au cas o\`u les pentes sont enti\`eres.
Cette op\'eration revient \`a une extension de corps aux
$q$-diff\'erences, soit encore \`a une extension du corps de base
de l'alg\`ebre $\mathcal{D}_{q}$. Elle est donc compatible avec
les op\'erations de $\mathcal{D}_{q}$. \\

\textsl{1.1.5 Changement de fonction inconnue $f = u g$. -}
Soit $u$ un \'el\'ement inversible de $L$ tel que 
$\sigma_{q}u = \alpha u \;,\; \alpha \in K^{*}$.
On a alors $P f = 0 \Leftrightarrow P^{[u]} g = 0$, avec:
$$
P^{[u]} \underset{def}{=}  u^{-1} P u = 
\sum_{i=0}^{n} a_{i} \frac{\sigma_{q}^{n-i}(u)}{u} \sigma^{n-i} =
\left(\prod_{i=0}^{n-1} \sigma_{q}^{i}(\alpha)\right)
\sum_{i=0}^{n} 
  \frac{a_{i}}{\sigma_{q}^{i}(\alpha)\cdots \sigma_{q}^{n-1}(\alpha)} 
  \sigma^{n-i}.
$$
Les pentes de $P^{[u]}$ sont donc
$\mu_{1} - v_{0}(\alpha),\ldots,\mu_{k} - v_{0}(\alpha)$.
On peut ainsi ramener une pente enti\`ere \`a $0$: 
si $\mu_{i} \in \mathbf{Z}$, on prend $u = \Theta_{q}^{\mu_{i}}$
(voir les conventions g\'en\'erales, dans l'introduction). \\

\textsl{1.1.6 Symboles de transformation de jauge. -} 
On peut remarquer que, d'apr\`es la troisi\`eme formule de 1.1.5,
$P^{[u]}$ peut \^etre d\'efini en fonction de $\alpha$ seul. 
On prouve alors directement:

\begin{itemize}

\item{Que $P \mapsto P^{[u]}$ est un automorphisme 
de $\mathcal{D}_{q}$. Il est en effet clair qu'il
est $\mathbf{C}$-lin\'eaire; la multiplicativit\'e
peut donc \^etre v\'erifi\'ee sur le produit
$(a \sigma^{i}).(b \sigma^{j}) = a \sigma_{q}^{i}(b) \sigma^{i+j}.$
Notons
$\pi_{i}(\alpha) = \frac{\sigma_{q}^{i}(u)}{u} =
\alpha \sigma_{q}(\alpha) \cdots \sigma_{q}^{i-1}(\alpha)$, 
de sorte que:
\begin{eqnarray*}
\left((a \sigma^{i}).(b \sigma^{j})\right)^{[u]} & = &
a \sigma_{q}^{i}(b) \pi_{i+j}(\alpha) \sigma^{i+j}           \\
(a \sigma^{i})^{[u]} (b \sigma^{j})^{[u]}        & = &
(a \pi_{i}(\alpha) \sigma^{i}).(b \pi_{j}(\alpha) \sigma^{j}) \\
                                                 & = &
a \sigma_{q}^{i}(b) 
\pi_{i}(\alpha) \sigma_{q}^{i}(\pi_{j}(\alpha)) \sigma^{i+j}
\end{eqnarray*}
Il s'agit donc de v\'erifier que 
$\pi_{i}(\alpha) \sigma_{q}^{i}(\pi_{j}(\alpha)) = \pi_{i+j}(\alpha)$,
ce qui est imm\'ediat (c'est un cocycle).}

\item{Que $P^{[uv]} = \left(P^{[u]}\right)^{[v]}$. 
Cela se ram\`ene, par des calculs encore plus simples, 
\`a v\'erifier que 
$\pi_{i}(\alpha \beta) = \pi_{i}(\alpha) \pi_{i}(\beta)$,
ce qui est \'evident.}

\end{itemize}

Nous garderons cependant la notation avec le ``symbole'' $u$, 
plus suggestive car elle rappelle qu'il s'agit en fait
d'une transformation de jauge. La seule r\`egle est 
que ces transformations agissent comme des automorphismes
int\'erieurs. \\

\textsl{1.1.7 Choix canonique des symboles. -}
Tout \'el\'ement $\alpha$ de $K^{*}$ s'\'ecrit de mani\`ere unique
$\alpha = c z^{\mu} \beta$, o\`u 
$c \in \mathbf{C}^{*}$, $\mu \in \mathbf{Z}$ et $\beta(0) = 1$.
Pour r\'esoudre l'\'equation $\sigma_{q} u = \alpha u$,
nous prendrons $u = e_{q,c} \Theta_{q}^{\mu} v$, o\`u
$v(z) = \underset{k \geq 1}{\prod} \beta(q^{-k}z)$;
on v\'erifie en effet facilement que $v \in \mathbf{K}^{*}$
(aussi bien dans le cas formel que dans le cas convergent).
On notera $e_{q,\alpha}$ l'\'el\'ement $u$ ainsi obtenu; 
c'est un \'el\'ement inversible de $L$.

\subsubsection*{Equation caract\'eristique, exposants}

On suppose ici que $S(P) \subset \mathbf{Z}$. On va d\'efinir
l'\'equation caract\'eristique et les exposants attach\'es
\`a la $i$-\`eme pente $\mu = \mu_{i}$ de $P$. 
D'apr\`es 1.1.5, la $i$-\`eme pente de 
$P^{[e_{q,z^{\mu}}]} = a'_{0} \sigma^{n} + \cdots + a'_{n}$
vaut $0$. Il existe donc des indices $\alpha \lt \beta$ 
dans $\{0,\ldots,n\}$ tels que,
notant $l = v_{0}(P^{[e_{q,z^{\mu}}]})$:
$$
\begin{cases}
v_{0}(a'_{\alpha}) = v_{0}(a'_{\beta}) = l \\
\forall i \in \{0,\ldots,n\} \;,\; v_{0}(a'_{i}) \geq l \\
\forall i \in \{0,\ldots,\alpha-1\} \cup \{\beta+1,\ldots,n\}
        \;,\; v_{0}(a'_{i}) \gt l
\end{cases}
$$
On a donc $r_{i} = \beta - \alpha$. On introduit: 
$Q = z^{-l} P^{[e_{q,z^{\mu}}]} = b_{0} \sigma^{n} + \cdots + b_{n}$,
dont les coefficients sont donc dans l'anneau de valuation $\mathcal{O}$
de $K$. Plus pr\'ecis\'ement, en posant:
\begin{eqnarray*} 
\overline{Q} & \underset{def}{=} & 
                                 b_{0}(0) \; \sigma^{n} + \cdots + b_{n}(0) \\
             & = & b_{\alpha}(0) \; \sigma^{n-\alpha} + \cdots + 
                   b_{\beta}(0) \; \sigma^{n-\beta} \\
             & \in & \mathbf{C}[\sigma,\sigma^{-1}]
\end{eqnarray*}
(qui est donc un polyn\^ome commutatif), on a 
$Q \equiv \overline{Q} \pmod{z \mathcal{O}\left<\sigma,\sigma^{-1}\right>}$
et $b_{\alpha}(0) b_{\beta}(0) \not= 0$.
L'\'equation $\overline{Q} = 0$ (ainsi que le polyn\^ome $\overline{Q}$ 
lui-m\^eme) est appel\'ee \emph{\'equation caract\'eristique} 
attach\'ee \`a la pente $\mu$ de $P$; on peut la consid\'erer comme 
d\'efinie \`a un facteur $\alpha \sigma^{k}$ pr\`es, 
$\alpha \in \mathbf{C}^{*}, k \in \mathbf{Z}$. On la notera
$\overline{P}^{(\mu)}$, ou simplement $\overline{P}$
dans le cas de la pente $\mu = 0$. Si $\mu \not\in S(P)$, 
l'\'equation caract\'eristique est une constante non nulle. 
En g\'en\'eral:
$$
\overline{P}^{(\mu)} = 
\left(z^{-v_{0}(P^{[e_{q,z^{\mu}}]})} 
P^{[e_{q,z^{\mu}}]}\right)_{z = 0}.
$$
On voit donc, dans tous les cas, que $r_{P}(\mu)$ est
\'egal au degr\'e absolu $\deg(\overline{P}^{(\mu)})$
de l'\'equation caract\'eristique. \\

\textsl{1.1.8 Lemme. -}
\emph{L'\'equation caract\'eristique est multiplicative: 
$$
\forall P_{1}, P_{2} \in \mathcal{D}_{q}
\;,\; \forall \mu \in \mathbf{Q} \;,\;
\overline{P_{1} P_{2}}^{(\mu)} = \overline{P_{1}}^{(\mu)} \overline{P_{2}}^{(\mu)}.
$$}
\Pr
En effet, c'est une cons\'equence imm\'ediate des propri\'et\'es
de la valuation et du degr\'e absolu dans $\mathcal{D}_{q}$
et de 1.1.6. 
\hfill $\Box$ \\

\textsl{1.1.9 Th\'eor\`eme. -}
\emph{Le polygone de Newton est additif. Pr\'ecis\'ement, 
$P_{1}$ et $P_{2}$ \'etant des op\'erateurs aux $q$-diff\'erences 
comme ci-dessus: 
\begin{eqnarray*}
r_{P_{1}P_{2}} & = & r_{P_{1}} + r_{P_{2}}, \\
N(P_{1} P_{2}) & = & N(P_{1}) + N(P_{2}).
\end{eqnarray*}}
\Pr
On d\'eduit la premi\`ere \'egalit\'e du lemme 1.1.8 et du fait que
$r_{P}(\mu) = \deg(\overline{P}^{(\mu)})$; la deuxi\`eme \'egalit\'e est
alors cons\'equence du lemme 1.1.3.
\hfill $\Box$ \\

Les racines non nulles de l'\'equation caract\'eristique 
attach\'ee \`a la pente $\mu$ de $P$ sont appel\'es les 
\emph{exposants} attach\'es \`a cette pente. Tout complexe
non nul s'\'ecrit de mani\`ere unique:
$$
c = q^{\epsilon(c)} \overline{c} \quad \text{ avec} \quad
\epsilon(c) \in \mathbf{Z} \quad \text{ et } 
\quad 1 \leq |\overline{c}| \lt |q|.
$$
Nous identifierons l'ensemble quotient 
$\mathbf{E}_{q} = \mathbf{C}^{*}/q^{\mathbf{Z}}$
(qui est une courbe elliptique) \`a la \emph{couronne fondamentale}
$\{z \in \mathbf{C}^{*} \;/\; 1 \leq |z| \lt |q|\}$,
qui en est un syst\`eme de repr\'esentants dans $\mathbf{C}^{*}$,
et le repr\'esentant $\overline{c}$ \`a la classe de $c$
modulo $q^{\mathbf{Z}}$. Nous dirons qu'un exposant $c$ 
attach\'e \`a la pente $\mu$ est \emph{non r\'esonnant} 
si $\epsilon(c)$ est maximal pour sa classe de congruence,
autrement dit, si aucun $c q^{l}$, $l \in \mathbf{N}^{*}$,
n'est un exposant (attach\'e \`a cette pente). \\

\textsl{1.1.10 Manipulations \'el\'ementaires 
sur l'\'equation caract\'eristique et sur les exposants. -}
Les propri\'et\'es suivantes d\'ecoulent imm\'ediatement de 1.1.4,
1.1.5 et 1.1.6.
 
\begin{enumerate}

\item{Si $\sigma_{q} u = z^{l} u$, l'\'equation 
caract\'eristique attach\'ee \`a la pente $\mu$ 
de $P$ est \'egale \`a l'\'equation caract\'eristique 
attach\'ee \`a la pente $\mu - l$ de $P^{[u]}$.}

\item{Si $\sigma_{q} u = c u \;,\; c \in \mathbf{C}^{*}$, 
et si l'on note $\overline{Q}(\sigma)$ l'\'equation 
caract\'eristique attach\'ee \`a la pente $\mu$ de $P$, 
l'\'equation caract\'eristique attach\'ee \`a la pente $\mu$ 
de $P^{[u]}$ est $\overline{Q}(c \sigma)$. 
On peut donc toujours ramener un exposant donn\'e \`a $1$ 
par une telle transformation de jauge..}

\item{Dans ce dernier cas, quitte \`a utiliser encore
une transformation de jauge avec $u = z^{l} \;,\; l \in \mathbf{Z}$, 
on peut m\^eme supposer $1$ non r\'esonnant.}

\end{enumerate}

\textsl{1.1.11 Remarque. -} La ramification ne change pas le calcul
des exposants. Cependant, ceux-ci ne sont intrins\`equement d\'efinis
qu'\`a un facteur de $q^{\mathbf{Z}}$ pr\`es, donc fournissent des
invariants dans la courbe elliptique $\mathbf{E}_{q}$, et celle-ci
est remplac\'ee, par ramification, par une courbe isog\`ene. Ce point
sera pr\'ecis\'e en 2.2.4.

% 1.2

\subsection{Factorisations formelle et convergente d'un 
op\'erateur aux $q$-diff\'erences}

En principe, la r\'esolution de l'\'equation (1) et la factorisation
de l'op\'erateur aux $q$-diff\'erences $P$ sont \'etroitement
imbriqu\'ees. Comme c'est surtout de la factorisation que nous avons
besoin, l'essentiel des r\'esultats de r\'esolution a \'et\'e report\'e
dans l'appendice.

\subsubsection*{Facteur droit associ\'e
               \`a un exposant non r\'esonnant}

D'apr\`es 1.1.10, on peut ramener, par l'interm\'ediaire
de transformations de jauge simples,
toute pente $\mu$ \`a $0$ et tout exposant $c$ attach\'e \`a 
cette pente \`a $1$. On peut m\^eme supposer que $1$ est non 
r\'esonnant, autrement dit, qu'aucun $q^{k} \;,\; k \in \mathbf{N}^{*}$
n'est un exposant attach\'e \`a la pente $0$. \\

\textsl{1.2.1 Lemme. -}
\emph{Supposons que $0$ est une pente de $P$ et que $1$ est un exposant 
non r\'esonnant attach\'e \`a cette pente. L'\'equation (1) admet alors 
une unique solution s\'erie formelle $f$ telle que $f(0) = 1$.} \\

\Pr
Comme pr\'ec\'edemment, on peut supposer $N(P)$ cal\'e de sorte
que la pente nulle soit sur l'axe des abcisses. Les coefficients $a_{i}$ 
admettent donc un d\'eveloppement en s\'erie:
$$
\forall i \in \{0,\ldots,n\} \;,\; 
a_{i} = \sum_{j \geq 0} a_{i,j} z^{j}.
$$
L'\'equation caract\'eristique attach\'ee \`a la pente $0$ 
est donc:  
$\overline{P} = a_{0,0} \; \sigma^{n} + \cdots + a_{n,0}$. 
On \'ecrit $f = \underset{m \geq 0}{\sum} f_{m} z^{m}$ 
l'inconnue. Alors:
$$
P.f = \sum_{l \geq 0} g_{l} z^{l},
$$
o\`u l'on a pos\'e:
$$
g_{l} = \left(\sum_{m = 0}^{l} F_{l-m}(q^{m})f_{m}\right)
\quad , \quad
F_{j}(X) = \sum_{i = 0}^{n} a_{n-i,j} X^{i}.
$$
Ce polyn\^ome est constitu\'e des coefficients qui contribuent
\`a la tranche d'ordonn\'ee $j$ dans l'int\'erieur du polygone de Newton.
Ainsi $F_{0} = \overline{P}$, d'o\`u $F_{0}(1) = 0$ et
$F_{0}(q^{m}) \not= 0$ pour $m \geq 1$ (puisque $1$ est
exposant non r\'esonnant). On trouve les coefficients par
r\'ecurrence en identifiant les $g_{l}$ \`a $0$. Comme $F_{0}(1) = 0$,
$f_{0}$ est arbitraire et, comme $F_{0}(q^{m}) \not= 0$ 
pour $m \geq 1$, les $f_{m}$ sont d\'etermin\'es inductivement
de mani\`ere unique. 
\hfill $\Box$ \\

\textsl{1.2.2 Lemme. -}
\emph{Supposons que $0$ est une pente de $P$. Soit $c$ un exposant 
non r\'esonnant de multiplicit\'e $m \geq 1$ attach\'e \`a la pente $0$. 
On a alors une factorisation:
$$
P = Q.(\sigma - c).u_{m}^{-1} \cdots (\sigma - c).u_{1}^{-1},
$$
o\`u $u_{1},\ldots,u_{m}$ sont des s\'eries formelles telles
que $u_{i}(0) = 1$. De plus, le polygone de Newton de $Q$ 
s'obtient en diminuant de $m$ la longueur de la pente 
horizontale de celui de $P$: $r_{Q} = r_{P} - m \delta_{0}$,
et les \'equations caract\'eristiques correspondantes 
v\'erifient: $\overline{P} = (X - c)^{m} \overline{Q}$
} \\

\Pr
Traitons d'abord le cas o\`u $c = 1$. 
Soit $u_{1} = f$, la s\'erie formelle obtenue au lemme 1.2.1. 
L'op\'erateur aux $q$-diff\'erences $P$ admet une factorisation:
$$
P = P_{1}.(\sigma - 1).u_{1}^{-1},
$$
o\`u l'on a pos\'e:
$$
P_{1} = \sum_{j=0}^{n-1} 
        \left(- \sum_{i=0}^{j} a_{n-i} \; \sigma_{q}^{i}(u_{1}) \right)
        \sigma^{j}.
$$
Ceci se v\'erifie par le calcul suivant, valable sans hypoth\`ese 
sur $u_{1}$:
\begin{eqnarray*}
P_{1} (\sigma - 1) 
  & = & - \sum_{0 \leq i \leq j \leq n-1} 
           a_{n-i} \; \sigma_{q}^{i}(u_{1})(\sigma^{j+1} - \sigma^{j}) \\
  & = & \sum_{j=0}^{n-1} \sum_{i=0}^{j} 
         a_{n-i} \; \sigma_{q}^{i}(u_{1})\sigma^{j} -
        \sum_{j=1}^{n} \sum_{i=0}^{j-1} 
         a_{n-i} \; \sigma_{q}^{i}(u_{1})\sigma^{j} \\
  & = & a_{n} u_{1} + 
        \sum_{j=1}^{n-1} 
         a_{n-j} \; \sigma_{q}^{j}(u_{1})\sigma^{j} -
        \sum_{i=0}^{n-1} a_{n-i} \; \sigma_{q}^{i}(u_{1})\sigma^{n} \\
  & = & \left(\sum_{j=0}^{n} a_{n-j} \; \sigma^{j}\right)u_{1} -
        \sum_{i=0}^{n} a_{n-i} \; \sigma_{q}^{i}(u_{1})\sigma^{n} \\
  & = & P u_{1} - \left(P(\sigma_{q}) u_{1}\right) \sigma^{n}.
\end{eqnarray*}

Le polygone de Newton de $P_{1}$ s'obtient en diminuant de $1$ 
la longueur de la pente horizontale de celui de $P$:
$r_{Q} = r_{P} - \delta_{0}$,
et les \'equations caract\'eristiques correspondantes 
v\'erifient:
$\overline{P} = (X - 1) \overline{P_{1}}$.
Il suffit alors d'it\'erer le processus pour obtenir 
la factorisation:
$$
P = Q.(\sigma - 1).u_{m}^{-1} \cdots (\sigma - 1).u_{1}^{-1},
$$
o\`u $u_{1},\ldots,u_{m}$ sont des s\'eries formelles telles
que $u_{i}(0) = 1$. De plus, le polygone de Newton de $Q$ 
s'obtient en diminuant de $m$ la longueur de la pente horizontale 
de celui de $P$, autrement dit, $r_{Q} = r_{P} - m \delta_{0}$. 
Les \'equations caract\'eristiques correspondantes v\'erifient: 
$\overline{P} = (X - 1)^{m} \overline{Q}$. \\

Dans le cas d'un exposant $c$ quelconque (non r\'esonnant),
soit $u = e_{q,c}$ (1.1.7). Alors $P^{[u]}$ v\'erifie les hypoth\`eses 
du premier cas. On \'ecrit donc:
$$
P^{[u]} = 
R.(\sigma - 1).u_{m}^{-1} \cdots (\sigma - 1).u_{1}^{-1}.
$$
On applique \`a cette \'egalit\'e la transformation de
jauge de symbole $u^{-1}$, qui commute au produit, qui
n'affecte pas les fonctions $u_{i}$ et qui transforme
$\sigma - 1$ en $c^{-1} \sigma - 1$. 
On obtient alors la factorisation voulue en prenant 
$Q = c^{-m} R^{[u^{-1}]}$. Le reste suit.
\hfill $\Box$ \\

\textsl{1.2.3 Proposition. -}
\emph{Soit $\mu$ une pente enti\`ere de $P$. Soit $c$ 
un exposant non r\'esonnant de multiplicit\'e $m$ attach\'e 
\`a la pente $\mu$. On a alors une factorisation:
$$
P = Q.(z^{-\mu} \sigma - c).u_{m}^{-1} \cdots
      (z^{-\mu} \sigma - c).u_{1}^{-1},
$$
o\`u $u_{1},\ldots,u_{m}$ sont des s\'eries formelles telles
que $u_{i}(0) = 1$. De plus, le polygone de Newton de $Q$ 
s'obtient en diminuant de $m$ la longueur
de la pente de valeur $\mu$ de celui de $P$:
$r_{Q} = r_{P} - m \delta_{\mu}$,
et les \'equations caract\'eristiques correspondantes 
v\'erifient:
$\overline{P}^{(\mu)} = (X - c)^{m} \overline{Q}^{(\mu)}$.} \\

\Pr
Soit $u = e_{q,z^{\mu}}$ (1.1.7). Alors $P^{[u]}$ v\'erifie 
les hypoth\`eses du lemme 1.2.2. On \'ecrit donc:
$$
P^{[u]} = 
R.(\sigma - c).u_{m}^{-1} \cdots (\sigma - c).u_{1}^{-1}.
$$
On applique \`a cette \'egalit\'e la transformation de
jauge de symbole $u^{-1}$, et l'on invoque 1.1.6.
\hfill $\Box$

\subsubsection*{Factorisation formelle 
d'un op\'erateur aux $q$-diff\'erences}

Il y a divers \'enonc\'es possibles, en voici un: \\

\textsl{1.2.4 Proposition. -}
\emph{Soit $\mu$ une pente enti\`ere de $P$. Soient 
$c_{1},\ldots,c_{p}$ les exposants attach\'es \`a 
la pente $\mu$, et $m_{1},\ldots,m_{p}$ leurs multiplicit\'es
respectives. On suppose les $c_{i}$ index\'es de telle
sorte que, si 
$\frac{c_{j}}{c_{i}} = q^{l} \;,\; l \in \mathbf{N}^{*}$,
alors $i \lt j$ (les exposants les moins r\'esonnants
sont factoris\'es \`a droite les premiers). On a alors une 
factorisation $P = Q R$, o\`u $\mu \not\in S(Q)$ et o\`u
$R = R_{1} \cdots R_{p}$, avec:
$$
\forall i \in \{1,\ldots,p\} \;,\;
R_{i} = (z^{-\mu} \sigma - c_{i}).u_{i,m_{i}}^{-1} \cdots
        (z^{-\mu} \sigma - c_{i}).u_{i,1}^{-1},
$$
les $u_{i,j}$ \'etant des s\'eries formelles telles que
$u_{i,j}(0) = 1$. On, de plus, 
$r_{R} = (m_{1} + \cdots + m_{p}) \delta_{\mu}$
et 
$\overline{R}^{(\mu)} = 
\underset{i=1}{\overset{p}{\prod}} (X - c_{i})^{m_{i}}$.} \\

\Pr
Il suffit d'appliquer r\'ep\'etitivement 
la proposition 1.2.3.
\hfill $\Box$ \\

On peut donner des conditions plus souples. Pour la r\'esolution 
(formelle ou convergente), il sera au contraire commode d'\^etre 
plus rigide et de consid\'erer une classe d'exposants \`a la fois 
(modulo $q^{\mathbf{Z}}$). On obtient de la m\^eme mani\`ere: \\

\textsl{1.2.5 Proposition. -}
\emph{Soit $\mu$ une pente enti\`ere de $P$. Soient 
$c_{1},\ldots,c_{p}$ les exposants d'une m\^eme classe
modulo $q^{\mathbf{Z}}$ attach\'es \`a la pente $\mu$, 
et $m_{1},\ldots,m_{p}$ leurs multiplicit\'es respectives. 
On suppose les $c_{i}$ index\'es de telle sorte que
$\epsilon(c_{1}) \lt \cdots \lt \epsilon(c_{p})$ 
(les exposants les moins r\'esonnants sont factoris\'es 
les premiers). On a alors une unique factorisation
$P = Q R$, o\`u aucun exposant attach\'e \`a la pente
$\mu$ de $Q$ n'est congru aux $c_{i}$ modulo $q^{\mathbf{Z}}$
et o\`u $R = R_{1} \cdots R_{p}$, avec:
$$
\forall i \in \{1,\ldots,p\} \;,\;
R_{i} = (z^{-\mu} \sigma - c_{i}).u_{i,m_{i}}^{-1} \cdots
        (z^{-\mu} \sigma - c_{i}).u_{i,1}^{-1},
$$
les $u_{i,j}$ \'etant des s\'eries formelles telles que
$u_{i,j}(0) = 1$. On, de plus, 
$r_{R} = (m_{1} + \cdots + m_{p}) \delta_{\mu}$
et 
$\overline{R}^{(\mu)} = 
\underset{i=1}{\overset{p}{\prod}} (X - c_{i})^{m_{i}}$.} \\
\hfill $\Box$

\subsubsection*{Factorisation convergente 
d'un op\'erateur aux $q$-diff\'erences}

Les r\'esultats qui pr\'ec\`edent, ainsi que leur application
\`a la r\'esolution formelle (appendice) sont dus \`a Adams. 
Mais l'\'enonc\'e le plus caract\'eristique de la th\'eorie, 
le \emph{lemme d'Adams}, est l'existence de \emph{solutions convergentes
associ\'ees \`a la premi\`ere pente} (\cite{Adams1}, \cite{Adams2}). 
Il se prouve le plus ais\'ement via la factorisation analytique,
bien que celle-ci ait \'et\'e obtenue ult\'erieurement
(\emph{voir} \cite{Birkhoff3}). \\

Nous prenons ici pour corps de base $K = \mathcal{M}(\mathbf{C})$
ou $K = \mathbf{C}(\{z\})$. Nous dirons alors qu'une s\'erie
(resp. une factorisation) est \emph{convergente} si elle d\'efinit
un \'el\'ement de $K$ (resp. si tous les facteurs sont \`a coefficients
dans $K$). \\

\textsl{1.2.6 Lemme. -}
\emph{On reprend d'abord les hypoth\`eses du lemme 1.2.1, 
en supposant de plus que $0$ est la premi\`ere pente, i.e. 
la plus petite. La s\'erie formelle $f$ obtenue comme solution 
est alors convergente.} \\

\Pr
Nous commen\c{c}ons par le cas o\`u $K = \mathbf{C}(\{z\})$.
Avec les conventions de 1.2.1, la premi\`ere pente vaut $0$,
les suivantes sont $\gt 0$ et:
$$
\begin{cases}
\deg \overline{P} = \deg F_{0} = n \\
\forall j \in \{1,\ldots,n\} \;,\; \deg F_{j} \leq n
\end{cases}
$$
Ce qui suit est alors une application de la m\'ethode des
s\'eries majorantes. Des conditions sur les degr\'es et
du fait qu'aucun $F_{0}(q^{l}) \;,\, l \geq 1$ ne s'annule, 
on tire:
$$
\exists A \gt 0 \;:\; \forall l \in \mathbf{N}^{*} \;,\; 
        \left|F_{0}(q^{l})\right| \geq A |q|^{ln}.
$$
De la convergence des coefficients $a_{0},\ldots,a_{n}$, on tire:
$$
\exists B,C \gt 0 \;:\; 
    \forall i \in \{0,\ldots,n\}, \forall j \in \mathbf{N} \;,\;
       |a_{i,j}| \leq B C^{j},
$$
d'o\`u l'on d\'eduit (avec la condition sur les degr\'es):
$$
\forall j \in \mathbf{N}^{*}, \forall l \in \mathbf{N} \;,\;
    \left|F_{j}(q^{l})\right| \leq (n+1) B C^{j} |q|^{ln}.
$$
Il vient, pour $l \geq 1$:
$$
|f_{l}| \leq 
\frac{(n+1) B 
\left(
C |q|^{(l-1)n} |f_{l-1}| + \cdots + C^{l} |q|^{(l-l)n} |f_{0}|
\right)}
{A |q|^{ln}}.
$$
On pose $g_{l} = \frac{|q|^{ln} |f_{l}|}{C^{l}}$ 
et $D = \frac{(n+1)B}{A}$, et l'on a:
$$
\begin{cases}
g_{0} = 1 \\
\forall l \geq 1 \;,\; g_{l} \leq D(g_{0} + \cdots + g_{l-1})
\end{cases}
$$
On montre alors par r\'ecurrence que $g_{l} \leq (D+1)^{l}$,
d'o\`u:
$$
\forall l \geq 0 \;,\; 
|f_{l}| \leq \left(\frac{C(D+1)}{|q|^{n}}\right)^{l}.
$$
On a donc bien $f \in \mathbf{C}(\{z\})$. \\

Supposons maintenant que $K = \mathcal{M}(\mathbf{C})$.
Nous appliquons le principe, \'egalement caract\'eristique 
des \'equations aux $q$-diff\'erences, selon lequel
``l'\'equation fonctionnelle propage la m\'eromorphie''.
Nous r\'e\'ecrivons l'\'equation (1) sous la forme:
$$
f(z) = 
- \sum _{i=1}^{n} 
\frac{a_{i}(q^{-n}z)}{a_{0}(q^{-n}z)} f(q^{-i}z).
$$
Si $f$ est d\'efinie et m\'eromorphe dans un disque de
centre $0$ et de rayon $r \gt 0$, et qu'elle y v\'erifie
cette \'equation, la m\^eme formule permet de la prolonger
en une fonction m\'eromorphe sur le disque de centre $0$ 
et de rayon $|q|r$, qui y v\'erifie la m\^eme \'equation.
On obtient ainsi (puisque $|q| \gt 1$) un prolongement
\`a $\mathbf{C}$ tout entier. 
\hfill $\Box$ \\

\textsl{1.2.7 Lemme. -}
\emph{On reprend les hypoth\`eses de 1.2.2,
en supposant de plus que $\mu$ est la premi\`ere pente. 
Les factorisations obtenues sont convergentes.} \\

\Pr
En effet, seul le calcul de $P_{1}$ dans la premi\`ere
\'etape est non formel, et il est clair qu'il fournit
un polyn\^ome en $\sigma$ \`a coefficients convergents. 
\hfill $\Box$ \\

\textsl{1.2.8 Th\'eor\`eme (Birkhoff-Guenther). -}
\emph{On reprend les hypoth\`eses de 1.2.4 et 1.2.5, en supposant
de plus que $\mu$ est la premi\`ere pente. Les factorisations 
obtenues en 1.2.4 et 1.2.5 sont convergentes.} \\

Encore une fois, il suffit de rassembler les morceaux.
\hfill $\Box$ \\

\textsl{1.2.9 Remarque. -}
Le rayon de convergence garanti par la preuve de 1.2.6
est $\frac{|q|}{C(D+1)}$ : le num\'erateur d\'epend 
de $q$ seul, le d\'enominateur de l'\'equation seule.
Notons par ailleurs que cette preuve est le seul point
qui d\'epend de l'hypoth\`ese $|q| > 1$ (voir aussi 
la remarque 3.1.3).

%%%%%%%%%%%%%%%%%%%%%%%%%%%%%%%%%%%%%%%%%%%%%%%%%%%%%%%%%%%%%%%%%%%%%%%%%%%%%%
 
% 2

\section{Polygone de Newton d'un module aux $q$-diff\'erences}

% 2.1

\subsection{Equations, syst\`emes, modules}

\subsubsection*{Les objets}

Nous d\'ecrivons ici diff\'erents mod\`eles de l'\'equation aux 
$q$-diff\'erences scalaire d'ordre $n$ et leurs relations. 
Nous partons de l'\'equation (1), dans laquelle nous supposons
$a_{0} = 1$. Notons, pour $f \in L$:
$$
X_{f} = 
\begin{pmatrix} 
f \\ \sigma_{q} f \\ \vdots \\ \sigma_{q}^{n-1} f 
\end{pmatrix}.
$$
Alors $f$ est solution de l'\'equation (1) 
si et seulement si $X_{f}$ est solution de l'\'equation (ou 
syst\`eme) lin\'eaire de rang $n$: 
\begin{equation}
\sigma_{q} X = A X,
\end{equation}
o\`u l'inconnue $X$ est un vecteur et o\`u la matrice du 
syst\`eme est:
$$
A =
\begin{pmatrix}
0 & 1 & 0 & \ldots & 0 & 0 \\
0 & 0 & 1 & \ldots & 0 & 0 \\
\vdots & \vdots & \vdots & \vdots & \vdots & \vdots \\
0 & 0 & 0 & \ldots & 1 & 0 \\
0 & 0 & 0 & \ldots & 0 & 1 \\
- a_{n} & - a_{n-1} & - a_{n-2} & \ldots & - a_{2} & - a_{1} 
\end{pmatrix}
\in GL_{n}(K).
$$
R\'eciproquement, il d\'ecoulera indirectement du lemme du vecteur 
cyclique (cf. \emph{infra} 2.1.1 et 2.1.11) que tout syst\`eme (2) 
est \emph{\'equivalent} (au sens expliqu\'e plus loin) \`a un 
syst\`eme provenant d'une \'equation scalaire (1). \\

On peut \'ecrire le syst\`eme (2) sous forme d'\'equation
au point fixe: $\Phi_{A}(X) = X$, o\`u l'on a d\'efini 
$\Phi_{A}: X \mapsto A^{-1} \left(\sigma_{q} X \right)$.
L'op\'erateur $\Phi_{A}$ est \emph{$\sigma_{q}$-lin\'eaire}  
autrement dit, il est additif et 
$\forall \lambda ,X \;:\: 
\Phi_{A}(\lambda X) = \sigma_{q}(\lambda) \Phi_{A}(X)$
\footnote{Cette repr\'esentation est analogue \`a la 
repr\'esentation des solutions du syst\`eme diff\'erentiel 
$\frac{dX}{dz} = A X$ comme \'el\'ements du noyau de la connexion 
$\Delta : X \mapsto \frac{dX}{dz} - A X$. La semi-lin\'earit\'e
correspond \`a la formule de Leibnitz.}. 
Ceci conduit \`a d\'efinir un \emph{module aux $q$-diff\'erences} 
sur le corps aux $q$-diff\'erences $(K,\sigma_{q})$ comme un 
couple 
$(M,\Phi)$, o\`u $M$ est un $K$-espace vectoriel de dimension 
finie et $\Phi$ un automorphisme $\sigma_{q}$-lin\'eaire de $M$. 
Par le choix d'une base de $M$, tout tel module s'identifie \`a 
un module de la forme $(K^{n},\Phi_{A})$ avec $A \in GL_{n}(K)$. \\

De m\^eme qu'un couple $(V,\phi)$ form\'e d'un $K$-espace
vectoriel $V$ et d'un endomorphisme (resp. d'un automorphisme)
$K$-lin\'eaire $\phi$ de $V$ peut \^etre consid\'er\'e comme
un module sur l'anneau principal $K[X]$ (resp. $K[X,X^{-1}]$),
de m\^eme un module aux $q$-diff\'erences $(M,\Phi)$ peut \^etre
consid\'er\'e comme un module \emph{\`a gauche} sur l'anneau 
$\mathcal{D}_{q}$ (cf. les conventions g\'en\'erales): 
si $P = \sum a_{i} \sigma^{i} \in \mathcal{D}_{q}$
et $x \in M$, on pose $P.x = P(\Phi)(x) = \sum a_{i} \Phi^{i}(x)$.
Le fait que $\dim_{K}(M) < \infty$ \'equivaut au fait que le 
$\mathcal{D}_{q}$-module $M$ est de longueur finie.
Les modules aux $q$-diff\'erences sont donc exactement les
$\mathcal{D}_{q}$-modules \`a gauche de longueur finie. 
Par abus de langage, on appellera \emph{rang} 
d'un $\mathcal{D}_{q}$-module \`a gauche de longueur finie,
ou du module aux $q$-diff\'erences correspondant,
la dimension du $K$-espace vectoriel sous-jacent. \\

Soit $F \in GL_{n}(K)$.  Si l'on pose $X = F Y$ dans (2), 
on est conduit \`a une \'equation $\sigma_{q} Y = B Y$, avec
$B = \left(\sigma_{q} F\right)^{-1} A F$,
que l'on consid\`erera comme \'equivalente \`a (2). On notera
en cons\'equence $A \sim  \left(\sigma_{q} F\right)^{-1} A F$.
Cela \'equivaut \`a l'isomorphie des $\mathcal{D}_{q}$-modules
associ\'es (voir la description intrins\`eque des morphismes
en 2.1.3). A titre d'exemple, si l'on part de l'\'equation d'ordre $1$:
$\sigma_{q} f = a f \;,\; a \in K^{*}$, on obtient le module
$(K,a^{-1} \sigma_{q})$. C'est aussi le $\mathcal{D}_{q}$-module
monog\`ene $\mathcal{D}_{q}/\mathcal{D}_{q} P$, o\`u
$P = \sigma - a^{-1}$. On voit facilement que tout module
de rang $1$ est de cette forme. De plus, l'\'equation 
$\sigma_{q} f = a f$ est
\'equivalente \`a l'\'equation $\sigma_{q} f = b f$ si et
seulement si 
$\exists u \in K^{*} \;:\; 
\frac{b}{a} = \frac{\sigma_{q}(u)}{u}$. \\

\textsl{2.1.1 Lemme du vecteur cyclique. -}
\emph{Soit $(M,\Phi)$ un module aux $q$-diff\'erences de rang $n$
sur $K$. Alors il existe $x \in M$ tel que 
$(x,\Phi(x),\ldots,\Phi^{n-1}(x))$ est une base de $M$.} \\

\Pr
Une preuve analytique, due \`a Birkhoff, est donn\'ee dans 
\cite{JSAIF} et une preuve alg\'ebrique plus g\'en\'erale,
inspir\'ee de Katz, est donn\'ee dans \cite{LDVpreprint}.
\hfill $\Box$ \\

On dira que $(x,\Phi(x),\ldots,\Phi^{n-1}(x))$ est une 
\emph{base cyclique} et que $x$ est un \emph{vecteur cyclique}. 
On peut alors \'ecrire 
$\Phi^{n}(x) + a_{1} \Phi^{n-1}(x) + \cdots + a_{n} x = 0$;
le polyn\^ome non commutatif
$P = \sigma^{n} + a_{1} \sigma^{n-1} + \cdots + a_{n}
\in \mathcal{D}_{q}$
annule le g\'en\'erateur $x$ du $\mathcal{D}_{q}$-module 
\`a gauche $M$. De l'euclidianit\'e de $\mathcal{D}_{q}$
on tire que $\mathcal{D}_{q} P$ est m\^eme l'id\'eal
annulateur de $x$: \\

\textsl{2.1.2 Corollaire. -}
\emph{Tout module aux $q$-diff\'erences est isomorphe \`a
un module de la forme $\mathcal{D}_{q}/\mathcal{D}_{q} P$,
o\`u $P$ est entier unitaire, autrement dit, de la forme 
ci-dessus.}
\hfill $\Box$ \\

On verra en 2.1.11 comment en d\'eduire qu'un tel module provient
d'une \'equation d'ordre $n$; il faut cependant prendre garde
que ce n'est pas l'\'equation $P.f = 0$ mais l'\'equation
\emph{duale}, comme on le voit dans l'exemple ci-dessus
(\'equation d'ordre $1$). Ceci sera d\'etaill\'e en 2.1.11. 
Plus g\'en\'eralement, soit $y$ un \'el\'ement du module
aux $q$-diff\'erences $M$. Il existe un plus grand entier
$p$ tel que la famille $(y,\Phi(y),\ldots,\Phi^{p-1}(y))$
est libre sur $K$; le sous-espace vectoriel qu'elle
engendre est stable par $\Phi$ et $\Phi^{-1}$, il admet donc
une structure de module aux $q$-diff\'erences. La relation
lin\'eaire 
$\Phi^{p}(y) + b_{1} \Phi^{p-1}(y) + \cdots + b_{p} y = 0$
d\'etermine l'unique polyn\^ome entier unitaire de degr\'e
minimum $Q = \sigma^{p} + b_{1} \sigma^{p-1} + \cdots + b_{p}$ 
tel que $Q(\Phi)(y) =0$: on l'appellera le \emph{polyn\^ome
(annulateur) minimal} de $y$. L'id\'eal \`a gauche annulateur 
dans $\mathcal{D}_{q}$ de $y$ est $\mathcal{D}_{q} Q$ et le
module engendr\'e par $y$ est isomorphe \`a 
$\mathcal{D}_{q}/\mathcal{D}_{q} Q$. La suite exacte
correspondante sera examin\'ee plus loin (cf. 2.2.2).

\subsubsection*{Les morphismes}

\textsl{2.1.3 Morphismes de modules aux $q$-diff\'erences. -} 
Heuristiquement, un morphisme du syst\`eme $\sigma_{q} X = A X$ de 
rang $n$ 
vers le syst\`eme $\sigma_{q} Y = B Y$ de rang $p$ devrait \^etre 
une transformation de jauge $X \mapsto Y = F X$ qui envoie les
solutions du premier vers les solutions du second. On le d\'efinit
donc comme une matrice $F \in M_{p,n}(K)$ telle que
$\left(\sigma_{q} F \right) A = B F$. 
De mani\`ere plus intrins\`eque, un morphisme de
$(M,\Phi)$ vers $(N,\Psi)$ est une application $K$-lin\'eaire
$f: M \rightarrow N$ telle que $\Psi \circ f = f \circ \Phi$
(un choix de bases redonne la description pr\'ec\'edente). 
Il est clair que cela \'equivaut pr\'ecis\'ement \`a la
$\mathcal{D}_{q}$-lin\'earit\'e pour l'application
correspondante entre $\mathcal{D}_{q}$-modules \`a gauche. 
La composition et les morphismes identit\'es sont d\'efinis
de mani\`ere \'evidente. \\

\textsl{2.1.4 Proposition. -}
\emph{On obtient ainsi une cat\'egorie ab\'elienne
$DiffMod(K,\sigma_{q})$, dans laquelle tout objet $M$
s'ins\`ere dans une suite exacte:
$$
0 \rightarrow \mathcal{D}_{q} \rightarrow \mathcal{D}_{q}
\rightarrow M \rightarrow 0.
$$
}
\Pr
En effet, la cat\'egorie $DiffMod(K,\sigma_{q})$ est 
la sous-cat\'egorie pleine de la cat\'egorie $\mathcal{D}_{q}Mod$
des modules \`a gauche sur l'anneau $\mathcal{D}_{q}$
form\'ee des modules de longueur finie. Par ailleurs, 
si $M = \mathcal{D}_{q}/\mathcal{D}_{q} P$, 
la fl\`eche $\mathcal{D}_{q} \rightarrow \mathcal{D}_{q}$
dans la suite exacte ci-dessus
n'est autre que la multiplication \`a droite par $P$.
\hfill $\Box$ \\

\textsl{2.1.5 Objet unit\'e et foncteurs solutions. -}
L'\'equation triviale $\sigma_{q} f = f$ a pour mod\`ele
le module aux $q$-diff\'erences $(K,\sigma_{q})$, soit
le $\mathcal{D}_{q}$-module \`a gauche
$\mathcal{D}_{q}/\mathcal{D}_{q}(\sigma - 1)$. C'est
\emph{l'unit\'e} pour la structure tensorielle introduite 
en 2.1.6 et on le note (traditionnellement) 
$\mathbf{\underline{1}}$. On d\'efinit les deux foncteurs 
repr\'esent\'es par $\mathbf{\underline{1}}$:
\begin{eqnarray*}
\Gamma(M) & \underset{def}{=} & Hom(\mathbf{\underline{1}},M) \\
\Gamma^{\vee}(M) & \underset{def}{=} & Hom(M,\mathbf{\underline{1}})
\end{eqnarray*}
Prenons pour $M$ un module aux $q$-diff\'erences admettant les
deux descriptions: $(M,\Phi) = (K^{n},\Phi_{A})$ 
et $\mathcal{D}_{q}/\mathcal{D}_{q} P$ (notations du d\'ebut de 
2.1). On a alors les identifications naturelles:
\begin{eqnarray*}
\Gamma(M) & = & \{x \in M \;/\; \Phi(x) = x\} 
              = \{X \in K^{n} \;/\; \sigma_{q} X = A X\} \\
\Gamma^{\vee}(M) & = & \{f \in K \;/\; P.f = 0\}
\end{eqnarray*}
Ces deux foncteurs sont exacts \`a gauche. Le foncteur
$\Gamma$ est covariant, c'est le \emph{foncteur des cosolutions},
le foncteur $\Gamma^{\vee}$ est contravariant, c'est le
\emph{foncteur des solutions} (ceci pour respecter ce qui semble
\^etre devenu l'usage standard en th\'eorie des
$\mathcal{D}$-modules). On peut aussi consid\'erer $\Gamma$
comme un foncteur ``sections globales'' (cette intuition est 
renforc\'ee par l'interpr\'etation faisceautique donn\'ee dans
\cite{RSZ}).

\subsubsection*{Constructions tensorielles}

Les constructions qui suivent et les preuves de leurs
propri\'et\'es sont d\'etaill\'ees dans \cite{SVdP} et
\cite{JSGAL}. La terminologie
tannakienne est celle de \cite{DF} et de \cite{DM}. \\

\textsl{2.1.6 Produit tensoriel et $Hom$ interne. -}
Soient $(M,\Phi)$ et $(N,\Psi)$ deux modules aux 
$q$-diff\'erences. Il y a un unique automorphisme 
$\sigma_{q}$-lin\'eaire de $M \otimes_{K} N$ tel que
$x \otimes y \mapsto \Phi(x) \otimes \Psi(y)$.
Il en fait un module aux $q$-diff\'erences. De m\^eme, 
l'automorphisme $\sigma_{q}$-lin\'eaire de $Hom_{K}(M,N)$ 
d\'efini par $f \mapsto \Psi \circ f \circ \Phi^{-1}$ en fait 
un module aux $q$-diff\'erences. On a donc deux bifoncteurs,
le produit tensoriel et le ``$Hom$ interne'', que l'on notera
$\underline{Hom}$. \\

\textsl{2.1.7 Lemme. -}
\emph{Ces deux foncteurs sont exacts, et on a la propri\'et\'e
d'adjonction:
$$
\underline{Hom}(M \otimes N,P) = 
\underline{Hom}(M,\underline{Hom}(N,P)).
$$
De plus, $\mathbf{\underline{1}}$ est une unit\'e pour 
le produit tensoriel. Le foncteur de passage au dual:
$M^{\vee} = \underline{Hom}(M,\mathbf{\underline{1}})$
est exact et compatible avec le produit tensoriel.}
\hfill $\Box$ \\

On voit que le $\mathbf{C}$-espace vectoriel $Hom(M,N)$
est en fait form\'e des sections globales (ou cosolutions)
de $\underline{Hom}(M,N)$. En particulier, le $K$-espace 
vectoriel sous-jacent au dual de $(M,\Phi)$ est le dual du 
$K$-espace vectoriel $M$ et le dual de $(K^{n},\Phi_{A})$
est $(K^{n},\Phi_{B})$, o\`u $B$ est la contragr\'ediente
${}^{t}A^{-1}$ de $A$. \\

\textsl{2.1.8 Proposition. -}
\emph{La cat\'egorie $DiffMod(K,\sigma_{q})$ est une cat\'egorie
tensorielle rigide $\mathbf{C}$-lin\'eaire.}
\hfill $\Box$ \\

Notons qu'il n'est pas facile d'obtenir un foncteur fibre
d\'efini sur $\mathbf{C}$ (\emph{voir} \cite{SVdP} pour la
voie alg\'ebrique, \cite{JSGAL} pour la voie analytique dans
le cas fuchsien).

\subsubsection*{Dualit\'e}

Les calculs qui suivent sont inspir\'es du lemme de Gabber
(cf. \cite{Katz2}, I.1.5). 
Soit $(M,\Phi)$ un module aux $q$-diff\'erences admettant la 
base cyclique $(e_{0},\ldots,e_{n-1})$. On a donc, pour 
$0 \leq i \leq n-2$, $\Phi(e_{i}) = e_{i+1}$ et
$\Phi(e_{n-1}) = - a_{n} e_{0} - \cdots - a_{1} e_{n-1}$,
o\`u $P = \sigma^{n} + a_{1} \sigma^{n-1} + \cdots + a_{n}$
est le polyn\^ome minimal de $e_{0}$. 
D'apr\`es 2.1.6 et 2.1.7, le dual de $(M,\Phi)$ est 
$(M^{\vee},\Phi^{\vee})$, o\`u $M^{\vee}$ est le dual de $M$ 
et $\Phi^{\vee}$ la contragr\'ediente ${}^{t}\Phi^{-1}$
de $\Phi$. Soit $(e_{0}^{\vee},\ldots,e_{n-1}^{\vee})$
la base duale de la base $(e_{0},\ldots,e_{n-1})$. On notera 
$\left<u,v\right>$ l'application d'une forme lin\'eaire $u$ 
\`a un \'el\'ement $v$.\\

\textsl{2.1.9 Lemme. -}
\emph{Soient $i,j \in \{0,\ldots,n-1\}$. Alors:
\begin{eqnarray*}
\text{(i) } 
\left<\left(\Phi^{\vee}\right)^{-1}(e_{i}^{\vee}),e_{j}\right>
& = &
\begin{cases}
\text{ si  } j \leq n-2 \;:\; \delta_{i,j+1} \\
\text{ si  } j = n-1 \;:\; - \sigma_{q}^{-1}(a_{n-i}) 
\end{cases} \\
\text{(ii) \hspace{0.9cm}}
\left(\Phi^{\vee}\right)^{-1}(e_{i}^{\vee})
& = &
\begin{cases}
\text{ si  } i = 0 \;:\; - \sigma_{q}^{-1}(a_{n}) e_{n-1}^{\vee} \\ 
\text{ si  } i \geq 1 \;:\;
e_{i-1}^{\vee} - \sigma_{q}^{-1}(a_{n-i}) e_{n-1}^{\vee}
\end{cases} \\
\text{(iii) \hspace{2.3cm}}
e_{i}^{\vee} & = & 
\begin{cases}
\text{ si  } i = 0 \;:\; - a_{n} \Phi^{\vee}(e_{n-1}^{\vee}) \\
\text{ si  } i \ \geq 1 \;:\; 
\Phi^{\vee}(e_{i-1}^{\vee}) - a_{n-i} \Phi^{\vee}(e_{n-1}^{\vee}) 
\end{cases}
\end{eqnarray*}
}
\Pr
La premi\`ere assertion se prouve en remarquant que, par 
d\'efinition du dual, on a la formule g\'en\'erale:
$\left<\Phi^{\vee}(u),\Phi(v)\right> = 
\sigma_{q}(\left<u,v\right>)$, d'o\`u
$\left<u,v\right> = 
\sigma_{q}^{-1}(\left<\Phi^{\vee}(u),\Phi(v)\right>)$.
On applique cette derni\`ere \'egalit\'e \`a
$u = \left(\Phi^{\vee}\right)^{-1}(e_{i}^{\vee})$ et 
$v = e_{j}$, ce qui donne:
\begin{eqnarray*}
\left<\left(\Phi^{\vee}\right)^{-1}(e_{i}^{\vee}),e_{j}\right> 
& = &
\sigma_{q}^{-1}(\left<e_{i}^{\vee},\Phi(e_{j})\right>) \\
& = &
\sigma_{q}^{-1}\left(
\begin{cases}
\text{ si  } j \leq n-2 \;:\; 
\left<e_{i}^{\vee},e_{j+1}\right> = \delta_{i,j+1} \\
\text{ si  } j = n-1 \;:\; 
\left<e_{i}^{\vee},- a_{n} e_{0} - \cdots - a_{1} e_{n-1}\right> 
= - a_{n-i} \\
\end{cases}
\right)
\end{eqnarray*}
d'o\`u l'on tire bien la formule annonc\'ee. \\

Pour la deuxi\`eme assertion, on part de la formule g\'en\'erale
$\forall u \in M^{\vee} \;,\;
u = \underset{j = 0}{\overset{n-1}{\sum}} \left<u,e_{j}\right>
e_{j}^{\vee}$, que l'on applique \`a 
$u = \left(\Phi^{\vee}\right)^{-1}(e_{i}^{\vee})$:
\begin{eqnarray*}
\left(\Phi^{\vee}\right)^{-1}(e_{i}^{\vee}) 
& = & 
\sum_{j=0}^{n-1} 
\left<\left(\Phi^{\vee}\right)^{-1}(e_{i}^{\vee}),e_{j}\right>
e_{j}^{\vee} \\
& = &
- \sigma_{q}^{-1}(a_{n-i}) e_{n-1}^{\vee} +
\begin{cases}
\text{ si  } i = 0    \;:\: 0 \\
\text{ si  } i \geq 1 \;:\; e_{i-1}^{\vee}
\end{cases}
\end{eqnarray*}
La troisi\`eme assertion vient alors imm\'ediatement.
\hfill $\Box$ \\

\textsl{2.1.10 Proposition. -}
\emph{Sous les hypoth\`eses ci-dessus, $e_{n-1}^{\vee}$
est un vecteur cyclique de $M^{\vee}$, de polyn\^ome minimal
$P^{\vee} = \sigma^{n} + b_{1} \sigma^{n-1} + \cdots + b_{n}$,
o\`u (en posant, par commodit\'e, $a_{0} = 1$):
$$
\forall i \in \{1,\ldots,n\} \;,\; 
b_{i} = 
\frac{\sigma_{q}^{n-i-1}(a_{n-i})}{\sigma_{q}^{n-1}(a_{n})}.
$$
}
\Pr
D'apr\`es le point (iii) du lemme 2.1.9, le plus petit
sous-espace de $M^{\vee}$ stable par $\Phi^{\vee}$ et contenant 
$e_{n-1}^{\vee}$ est $M^{\vee}$, donc $e_{n-1}^{\vee}$ est un 
vecteur cyclique, dont il reste \`a d\'eterminer le polyn\^ome
minimal. On d\'eduit \'egalement du lemme, par r\'ecurrence:
$$
\forall i \in \{0,\ldots,n-1\} \;:\;
e_{i}^{\vee} = 
- \sum_{j=0}^{i} \sigma_{q}^{i-j}(a_{n-j}) 
\left(\Phi^{\vee}\right)^{i-j+1}(e_{n-1}^{\vee}),
$$
d'o\`u $Q(\sigma_{q})(e_{n-1}^{\vee}) = 0$, avec
$Q = 1 + \underset{j = 0}{\overset{n - 1}{\sum}}
\sigma_{q}^{i-j}(a_{n-j}) \sigma^{n-j} \in \mathcal{D}_{q}$.
En divisant ce polyn\^ome (entier) par son coefficient
dominant $\sigma_{q}^{n-1}(a_{n})$, on obtient le polyn\^ome
minimal entier unitaire $P^{\vee}$ annonc\'e.
\hfill $\Box$ \\

\textsl{2.1.11 Lien entre solutions et cosolutions.}
On va maintenant \'elucider le lien entre les deux mod\`eles
de l'\'equation (1). On peut supposer celle-ci \'ecrite
sous la forme $P.f = 0$, o\`u 
$P = \sigma^{n} + a_{1} \sigma^{n-1} + \cdots + a_{n}$. 
Notons ${}_{P}M$ le module obtenu par transformation de (1)
en syst\`eme. C'est donc $K^{n}$ muni de l'automorphisme
$\sigma_{q}$-lin\'eaire 
$X \mapsto A^{-1} \left(\sigma_{q} X\right)$, o\`u $A$ est la
matrice d\'ecrite au d\'ebut de  2.1. 
Notons d'autre part $M_{P}$ le module 
$\mathcal{D}_{q}/\mathcal{D}_{q} P$. Son op\'erateur $\Phi$
est la multiplication par $\sigma$ modulo $P$. Une base
cyclique est donc form\'ee des classes modulo $P$:
$\overline{1},\overline{\sigma},\ldots,\overline{\sigma^{n-1}}$.
La matrice de $\Phi$ dans cette base est ${}^{t}A$. Ce module
est donc isomorphe \`a $K^{n}$ muni de l'automorphisme
$\sigma_{q}$-lin\'eaire 
$X \mapsto B^{-1} \left(\sigma_{q} X\right)$, o\`u
$B = {}^{t}A^{-1}$ est la contragr\'ediente de $A$. 
Ainsi, les mod\`eles ${}_{P}M$ et $M_{P}$ de (1) sont duaux
l'un de l'autre (2.1.7), ce qui explique pourquoi les solutions de
(1) peuvent \^etre vues au choix comme solutions de l'un ou
cosolutions de l'autre. On en d\'eduit aussi que tout module
aux $q$-diff\'erences est isomorphe \`a un module ${}_{P}M$,
autrement dit qu'il provient d'une \'equation.

% 2.2

\subsection{D\'efinition intrins\`eque du polygone de Newton}

\subsubsection*{D\'evissage et triangularisation}

Voici quelques premi\`eres pr\'ecisions sur la structure des modules
aux $q$-diff\'erences dans les cas formel et convergent. Seules les
propri\'et\'es n\'ecessaires \`a la d\'efinition du polygone de
Newton sont indiqu\'ees ici (voir cependant la remarque 2.2.4). \\

Notons $K^{*}_{\sigma_{q}}$ le groupe quotient de $K^{*}$ par
le sous-groupe $\{\frac{\sigma_{q}(u)}{u} \;/\; u \in K^{*}\}$. \\

\textsl{2.2.1 Lemme: classes de modules de rang $1$. -}
\emph{L'application qui, \`a $a \in K^{*}$ associe le module
$M_{\sigma - a}$ induit une bijection de $K^{*}_{\sigma_{q}}$ sur
l'ensemble des classes d'isomorphie de modules aux $q$-diff\'erences 
de rang $1$ sur $K$ .} \\

\Pr
Il d\'ecoule du lemme du vecteur cyclique (2.1.1) 
que tout objet $M$ de $DiffMod(K,\sigma_{q})$ s'\'ecrit sous la 
forme $M_{P} = \mathcal{D}_{q}/\mathcal{D}_{q} P$ avec $P$ entier 
unitaire (l'anneau $\mathcal{D}_{q}$ n'\'etant pas commutatif,
on ne peut cependant pas caract\'eriser l'id\'eal 
$\mathcal{D}_{q} P$ comme id\'eal annulateur de $M$ et en 
d\'eduire l'unicit\'e de $P$).
Par ailleurs, le rang du module $M$ (c'est \`a dire sa dimension
sur $K$, cf. le d\'ebut de 2.1) est un invariant d'isomorphie (car toute
bijection $\mathcal{D}_{q}$-lin\'eaire est $K$-lin\'eaire) et
$M_{P}$ a pour rang le degr\'e de $P$. Donc $M_{P}$ est de rang
$1$ si et seulement si $P = \sigma - a, a \in K^{*}$. De plus,
la description donn\'ee au d\'ebut de 2.1.1 entraine que tout
$P'$ tel que $M_{P} \simeq M_{P'}$ est, dans ce cas, de la forme
$P' = \sigma - a'$, o\`u
$\exists u \in K^{*} \;:\; 
\frac{a'}{a} = \frac{\sigma_{q}(u)}{u}$. 
\hfill $\Box$ \\

\textsl{2.2.2 Lemme: suites exactes et polyn\^omes minimaux. -}
\emph{Toute suite exacte dans $DiffMod(K,\sigma_{q})$ est isomorphe
\`a une suite exacte de la forme:
\begin{equation}
0 \rightarrow \mathcal{D}_{q}/\mathcal{D}_{q} Q \rightarrow
\mathcal{D}_{q}/\mathcal{D}_{q} P \rightarrow
\mathcal{D}_{q}/\mathcal{D}_{q} R \rightarrow 0,
\end{equation}
o\`u $P,Q,R$ sont des polyn\^omes non commutatifs entiers unitaires
tels que $P = QR$.} \\

\Pr
Soient $P,Q,R$ des polyn\^omes non commutatifs entiers unitaires
tels que $P = QR$. 
Alors $\mathcal{D}_{q} P \subset \mathcal{D}_{q} R$, d'o\`u
la surjection canonique 
$\mathcal{D}_{q}/\mathcal{D}_{q} P \rightarrow
\mathcal{D}_{q}/\mathcal{D}_{q} R \rightarrow 0$.
L'anneau $\mathcal{D}_{q}$ \'etant int\`egre, le noyau
$\mathcal{D}_{q} R/\mathcal{D}_{q} P =
\mathcal{D}_{q} R/\mathcal{D}_{q} QR$
est canoniquement isomorphe, via l'application lin\'eaire
\`a gauche $F \mapsto F R$, \`a 
$\mathcal{D}_{q}/\mathcal{D}_{q} Q$. On obtient donc une suite 
exacte de la forme indiqu\'ee.
R\'eciproquement, on d\'eduit facilement de 2.1 que toute
suite exacte dans $DiffMod(K,\sigma_{q})$ s'identifie \`a une 
suite exacte (3). 
\hfill $\Box$ \\

Tout objet $M$ de la cat\'egorie ab\'elienne 
$DiffMod(K,\sigma_{q})$ \'etant, par hypoth\`ese, de longueur
finie, il admet un d\'evissage:
$$
\{0\} = M_{0} \subset M_{1} \subset \cdots \subset M_{r} = M,
$$
o\`u tous les quotients $S_{i} = M_{i}/M_{i-1}$ sont simples.
D'apr\`es le th\'eor\`eme
de Jordan-H\"{o}lder, les classes d'isomorphie des modules $S_{i}$
sont bien d\'etermin\'ees dans leur ensemble, c'est \`a dire \`a
permutation pr\`es, et en conservant leurs multiplicit\'es. 
Le rang de $M$ est la somme des rangs des $S_{i}$.
Tout objet de rang $1$ est donc simple. 
Nous dirons que le module $M$ est \emph{triangularisable} si,
pour l'un des d\'evissages ci-dessus (et donc pour tous),
tous les quotients $S_{i} = M_{i}/M_{i-1}$ sont de rang $1$. 
Cela revient exactement \`a dire que la matrice
$A$ qui apparait dans la description sous la forme 
$(K^{n},\Phi_{A})$ peut \^etre choisie triangulaire sup\'erieure. \\

\textsl{2.2.3 Th\'eor\`eme. -}
\emph{Soit $M$ un module aux $q$-diff\'erences sur l'un des
corps aux $q$-diff\'erences $(K,\sigma_{q})$ mentionn\'es
dans les conventions g\'en\'erales. Il existe une extension
$(K_{l},\sigma_{q_{l}})$ de $(K,\sigma_{q})$ obtenue par 
ramification $z = z_{l}^{l}, q = q_{l}^{l}$ (cf. 1.1.4)
telle que le module $M_{l} = K_{l} \otimes_{K} M$ obtenu
par extension des scalaires est triangularisable.} \\

\Pr
En effet, c'est une cons\'equence imm\'ediate des r\'esultats
de 1.2 et du lemme 2.2.2 (comparer \`a \cite{Praagman}).
\hfill $\Box$ \\

\textsl{2.2.4 Remarque. -}
Le \emph{groupe} $K^{*}_{\sigma_{q}}$ est isomorphe au terme droit 
de la suite exacte:
$$
0 \rightarrow C_{K}^{*} \rightarrow K^{*} \rightarrow K^{*} 
\rightarrow K^{*}_{\sigma_{q}} \rightarrow 0,
$$
dans laquelle la fl\`eche centrale est l'application
$u \mapsto \frac{\sigma_{q}(u)}{u}$.
Les classes d'isomorphie de modules aux $q$-diff\'erences 
de rang $1$ sur $K$ forment un groupe pour le produit tensoriel,
et la bijection obtenue en 2.2.1 est un isomorphisme
(voir la preuve du th\'eor\`eme 2.3.1). 
Notons $C_{l}$ le groupe des classes de modules de rang $1$ sur
$(K_{l},\sigma_{q_{l}})$ (donc $C_{1} \simeq K^{*}_{\sigma_{q}}$)
et $C_{\infty}$ la limite inductive des $C_{l}$. 
A tout module aux $q$-diff\'erences $M$ on peut associer
une combinaison lin\'eaire formelle d'\'el\'ements de 
$C_{\infty}$ \`a coefficients dans $\mathbf{N}$; cette 
application est additive pour les suites exactes. 
On peut d\'eterminer pr\'ecis\'ement la structure de $C_{\infty}$.
Ses \'el\'ements portent \`a la fois l'information sur les pentes
et sur les exposants
\footnote{Dans le cas de pentes non enti\`eres, l'information
``galoisienne'' est cependant perdue: il vaut en fait mieux
consid\'erer le groupe $Div(C_{l})$
des \emph{combinaisons lin\'eaires formelles d'\'el\'ements 
de $C_{l}$ invariantes par $Gal(K_{l}/K)$}.}. 
La construction du polygone de Newton ne retient
que l'information sur les pentes.

\subsubsection*{Construction du polygone de Newton}

Notons $v_{0}$ la valuation $z$-adique sur les corps qui nous
int\'eressent. Notre but est d'attribuer la pente $v_{0}(a)$ 
\`a l'\'equation $\sigma_{q} f = a f$ et d'identifier celle-ci 
au module $\mathcal{D}_{q}/\mathcal{D}_{q} (\sigma - a)$, donc au
\emph{second} mod\`ele (celui des solutions et non celui des 
cosolutions)
\footnote{Ce point est important et la confusion entre 
les deux mod\`eles a entrain\'e la pr\'esence d'\'enonc\'es
inexacts dans ma note de C.R.A.S. sur ce sujet.}. \\

\textsl{2.2.5 L'algorithme. -}
La m\'ethode employ\'ee est inspir\'ee de Katz (cf \cite{Katz2}, 
II.2). On part d'un module aux $q$-diff\'erences $M$ sur 
$(K,\sigma_{q})$. On va construire la \emph{fonction de Newton} 
$r_{M} : \mathbf{Q} \rightarrow \mathbf{N}$ de $M$ (cf. 1.1.1).

\begin{enumerate}

\item{Quitte \`a ramifier, on peut supposer $M$ triangularisable.}

\item{Chaque module simple $S_{i}$ qui intervient dans la 
d\'ecomposition de $M$ est de rang $1$, donc de la forme
$M_{\sigma - a_{i}}$, l'\'el\'ement $a_{i}$ \'etant d\'etermin\'e \`a
un facteur $\frac{\sigma_{q}(u)}{u}$ pr\`es; en particulier,
la valuation $\mu_{i} = v_{0}(a_{i})$ est bien d\'etermin\'ee. On 
attribue \`a $S_{i}$ la fonction de Newton $r_{S_{i}} = \delta_{\mu_{i}}$:
autrement dit, le polygone de Newton de $S_{i}$ a une seule
pente, de valeur $\mu_{i}$ et de multiplicit\'e $1$.}

\item{On attribue \`a $M$ la fonction de Newton $r_{M}$ somme 
des $r_{S_{i}}$. D'apr\`es le th\'eor\`eme de Jordan-H\"{o}lder,
celle-ci est bien d\'etermin\'ee (i.e. ne d\'epend pas du 
d\'evissage choisi).}

\item{Si l'on a d\^u ramifier au niveau $l$, on divise toutes
les pentes calcul\'ees par $l$. Comme la ramification au niveau 
$l$ multiplie les valuations par $l$, le choix du niveau de 
ramification n'influe pas sur le r\'esultat final.}

\end{enumerate}

Cette d\'efinition du polygone de Newton est tautologiquement
additive pour les suites exactes (voir le lemme 1.1.3). 
Par construction, si $M_{P}$
est de rang $1$, son polygone de Newton est \'egal \`a celui du
polyn\^ome non commutatif $P$. D'apr\`es 1.1.9 et 2.2.2,
cela est encore vrai si $M_{P}$ est triangularisable. Enfin, 
d'apr\`es 1.1.2 et le dernier point de l'algorithme, cela reste
vrai en toute g\'en\'eralit\'e. \\

\textsl{2.2.6 Th\'eor\`eme et d\'efinition. -}
\emph{On peut associer \`a tout module aux $q$-diff\'erences
$M$ un polygone de Newton $N(M)$ d\'efini par la fonction
de Newton associ\'ee $r_{M}$, de mani\`ere que: \\
(i) Le polygone de Newton d'un polyn\^ome non commutatif entier
unitaire $P$ est celui de $M_{P}$: $N(M_{P}) = N(P)$. \\
(ii) Le passage au polygone de Newton est additif pour les
suites exactes.}
\hfill $\Box$ \\

\textsl{2.2.7 Exemple: premi\`ere pente d'un \'el\'ement. -}
Reprenant les notations qui suivent 2.1.2,
notons $Q = \sigma^{p} + b_{1} \sigma^{p-1} + \cdots + b_{p}$ 
le polyn\^ome annulateur minimal de $y \in M$. La premi\`ere pente de $Q$
est aussi la premi\`ere pente du sous-module $\mathcal{D}_{q} y$ de $M$.
Notons la $\mu(y)$. On v\'erifie facilement qu'elle vaut:
$$
\mu(y) = \underset{1 \leq i \leq n}{\min} \frac{v_{0}(b_{i})}{i}.
$$

% 2.3

\subsection{Propri\'et\'es fonctorielles, 
ab\'eliennes et tensorielles}

\subsubsection*{Propri\'et\'es du polygone de Newton}

Nous exprimerons surtout ces propri\'et\'es \`a l'aide de la
fonction de Newton. \\

\textsl{2.3.1 Th\'eor\`eme. -}
\emph{
(i) Le polygone de Newton est additif pour les suites exactes.
Pr\'ecis\'ement, si la suite de modules aux $q$-diff\'erences:
$$
0 \rightarrow M' \rightarrow M \rightarrow M'' \rightarrow 0
$$
est exacte, on a: 
$$
r_{M} = r_{M'} + r_{M''} \text{~et~} N(M) = N(M') + N(M'').
$$
(ii) Le polygone de Newton est multiplicatif par rapport au
produit tensoriel. Pr\'ecis\'ement,
si $M_{1}$ et $M_{2}$ sont des modules aux $q$-diff\'erences :
$$
\forall \mu \in \mathbf{Q} \;:\;
r_{M_{1} \otimes M_{2}}(\mu) =
\sum_{\mu_{1} + \mu_{2} = \mu} 
r_{M_{1}}(\mu_{1}) r_{M_{2}}(\mu_{2}).
$$
(iii) Le polygone de Newton du dual d'un module aux 
$q$-diff\'erences $M$ est sym\'etrique de celui de $M$;
pr\'ecis\'ement, sa fonction de Newton est donn\'ee par la formule:
$$
\forall \mu \in \mathbf{Q} \;:\; 
r_{M^{\vee}}(\mu) = r_{M}(-\mu).
$$
}
\Pr
La premi\`ere assertion a d\'ej\`a \'et\'e vue au 2.2. Pour la
seconde assertion (o\`u l'on a bien affaire \`a une somme finie 
de termes non nuls),  on v\'erifie d'abord la formule dans le cas 
o\`u $M_{1}$ et $M_{2}$ sont de rang $1$: elle vient alors, par 
application de la d\'efinition (2.2.5, point 2), de 
l'\'egalit\'e:
$$
\left(\mathcal{D}_{q}/\mathcal{D}_{q} (\sigma - a_{1})\right)
\otimes
\left(\mathcal{D}_{q}/\mathcal{D}_{q} (\sigma - a_{2})\right)
=
\mathcal{D}_{q}/\mathcal{D}_{q} (\sigma - a_{1}a_{2}),
$$
que l'on v\'erifie facilement. Le cas de deux modules 
triangularisables s'en d\'eduit gr\^ace \`a (i) et \`a 
l'exactitude du produit tensoriel (cf. 2.1.7).
Le cas g\'en\'eral vient alors de 2.2.3, du comportement
du polygone de Newton par ramification (cf. 1.1.2) et de
la compatibilit\'e du produit tensoriel avec l'extension des
scalaires (donc avec la ramification). Pour la troisi\`eme 
assertion, on raisonne de la m\^eme mani\`ere \`a partir de
l'\'egalit\'e:
$$
\left(\mathcal{D}_{q}/\mathcal{D}_{q} (\sigma - a)\right)^{\vee}
=
\mathcal{D}_{q}/\mathcal{D}_{q} (\sigma - a^{-1}).
$$
\hfill $\Box$ \\

\textsl{2.3.2 Remarque. -}
Ces propri\'et\'es s'expriment agr\'eablement \`a l'aide de la
s\'erie g\'en\'eratrice :
$$
\mathcal{R}_{M}(T) = \sum_{\mu \in \mathbf{Q}} r_{M}(\mu) T^{\mu},
$$
qui est un polyn\^ome ramifi\'e. On verra en 3.3 que l'on peut
en fait l'interpr\'eter comme la s\'erie de Hilbert-Samuel
d'un espace vectoriel gradu\'e par $\mathbf{Q}$. \\

Comme en 1.1.1, on introduit \emph{l'ensemble $S(M)$ des pentes
du module aux $q$-diff\'erences $M$}: c'est le support de sa
fonction de Newton $r_{M}$. Ce qui suit est imm\'ediat: \\

\textsl{2.3.3 Corollaire. -}
\emph{
(i) Si la suite de modules aux $q$-diff\'erences:
$$
0 \rightarrow M' \rightarrow M \rightarrow M'' \rightarrow 0
$$
est exacte, on a: 
$$
S(M) = S(M') \cup S(M'').
$$
(ii) Si $M_{1}$ et $M_{2}$ sont des modules aux $q$-diff\'erences:
$$
S(M_{1} \otimes M_{2}) = S(M_{1}) + S(M_{2}).
$$
(iii) L'ensemble des pentes du dual de $M$ est donn\'e par la 
formule: 
$$
S(M^{\vee}) = - S(M).
$$}
\hfill $\Box$ \\

Notons en particulier les tr\`es utiles cons\'equences 
suivantes: \\

\textsl{2.3.4 Corollaire. -}
\emph{
(i) Soit $N$ un quotient de $M$, par exemple, l'image $f(M)$
de $M$ par un morphisme. Alors $S(N) \subset S(M)$. \\
(ii) Soient $M'$ et $M''$ des sous-modules de $M$. Alors
$S(M') \;,\; S(M'') \subset S(M)$ et 
$S(M' + M'') = S(M') \cup S(M'')$. Si de plus 
$S(M') \cap S(M'') = \emptyset$, leur somme est directe. \\
(iii) Si $S(M) \cap S(N) = \emptyset$, tout morphisme de $M$
dans $N$ est nul.}
\hfill $\Box$

\subsubsection*{Modules purs, modules fuchsiens}

Un module \emph{pur de pente $\mu$} est un module $M$ tel que
$S(M) = \{\mu\}$. Notons que la question difficile, \`a ce stade, 
est de savoir si un module $M$ tel que $\mu \in S(M)$ admet un 
sous-module de pente $\mu$ (sous-entendu: non trivial), autrement 
dit, si l'on peut ``casser les pentes'' de son polygone de Newton.
Ce sera l'objet de la section 3.1.\\

\textsl{2.3.5 Proposition. -}
\emph{
(i) Tout sous-module, tout module quotient, toute image par un 
morphisme d'un module pur de pente $\mu$ sont soit triviaux, soit 
des modules purs de pente $\mu$. \\
(ii) Toute extension de modules purs de pente $\mu$ en est un. \\
(iii) Si $\mu \in S(M)$, la somme $N$ des sous-modules purs de 
pente $\mu$ de $M$ est soit nulle, soit le plus grand sous-module
pur de pente $\mu$ de $M$. Les sous-modules ainsi associ\'es aux
diff\'erentes pentes $\mu \in S(M)$ sont en somme directe. \\
(iv) Le dual d'un module pur de pente $\mu$ est un module pur de
pente $- \mu$. \\
(v) Le produit tensoriel de deux modules purs de pentes $\mu$ et
$\nu$ est un module pur de pente $\mu + \nu$.} \\

\Pr
Cela d\'ecoule imm\'ediatement de 2.3.3 et 2.3.4.
\hfill $\Box$ \\

Un module \emph{fuchsien} est un module pur de pente $0$. Le lien
avec les autres caract\'erisations des \'equations fuchsiennes et
des syst\`emes fuchsiens est explicit\'e dans l'appendice de
\cite{JSAIF} et dans \cite{JSGAL}. Leur r\^ole est crucial dans
l'\'etude g\'en\'erale des \'equations aux $q$-diff\'erences. \\

\textsl{2.3.6 Th\'eor\`eme. -}
\emph{La sous-cat\'egorie pleine de la cat\'egorie 
$DiffMod(K,\sigma_{q})$ form\'ee des modules fuchsiens est stable
par passage aux sous-quotients, aux extensions, au produit
tensoriel, au dual et aux $Hom$ internes. C'est donc une 
cat\'egorie tensorielle rigide $\mathbf{C}$-lin\'eaire
} \\

Il faut simplement noter que le $Hom$ interne s'exprime \`a l'aide
du produit tensoriel et du dual:
$$
\underline{Hom}(M,N) = M^{\vee} \otimes N.
$$
Tout le reste est cons\'equence imm\'ediate de 2.3.5. 
\hfill $\Box$ \\

Il est d\'emontr\'e dans \cite{SVdP} et dans \cite{JSGAL} 
que cette cat\'egorie est
tannakienne sur $\mathbf{C}$. Notons que, dans le th\'eor\`eme
ci-dessus, le plus compliqu\'e est de prouver qu'un sous-module
d'un module fuchsien est fuchsien (le cas du quotient s'en
d\'eduisant par dualit\'e). \\

\textsl{2.3.7 Formes canoniques. -}
Tout module pur de pente enti\`ere $\mu$ est le produit tensoriel
du module de rang $1$ (n\'ecessairement pur) et de pente $\mu$ :
$(K,z^{-\mu}\sigma_{q}) = 
\mathcal{D}_{q}/\mathcal{D}_{q} (\sigma - z^{\mu})$
par un module fuchsien. Des formes canoniques pour les modules
fuchsiens sont d\'ecrites dans \cite{JSAIF}. Le probl\`eme des 
formes canoniques pour les modules purs de pente non enti\`ere 
est int\'eressant et plus compliqu\'e. Il sera abord\'e dans 
\cite{JSIRR}.

%%%%%%%%%%%%%%%%%%%%%%%%%%%%%%%%%%%%%%%%%%%%%%%%%%%%%%%%%%%%%%%%%%%%%%%%%%%%%%

% 3

\section{La filtration canonique par les pentes
et le gradu\'e associ\'e}

% 3.1

\subsection{La filtration canonique par les pentes}

\subsubsection*{Le plus grand sous-module pur de pente $\mu$}

Soit $\mu \in S(M)$. On notera (uniquement dans ce paragraphe) 
$M^{[\mu]}$ le plus grand sous-module pur de pente $\mu$ de $M$, 
dont l'existence a \'et\'e \'etablie en 2.3.5,(iii). Il d\'ecoule
d'ailleurs de 2.3.5 et 2.3.1 qu'il est \emph{invariant par tout 
automorphisme du module aux $q$-diff\'erences $M$} d'une part, 
que son rang est \emph{a priori} major\'e par $r_{M}(\mu)$ d'autre part. \\

\textsl{3.1.1 Th\'eor\`eme. -}
\emph{Soit $\mu$ une pente du module $M$, suppos\'ee quelconque 
dans le cas formel, maximale ($\mu = \max S(M)$) dans le cas 
convergent. Alors $M$ admet un sous-module $M^{[\mu]}$ pur de 
pente $\mu$ et de rang maximum $r_{M}(\mu)$.} \\

\Pr
On suppose dans un premier temps que $S(M) \subset \mathbf{Z}$;
en particulier, $M$ est triangularisable (2.2.3). On \'ecrit
$M = M_{P}$ avec $P$ entier unitaire. On a donc 
$\mu \in S(M) = S(P)$. D'apr\`es 1.2.4 ou 1.2.5 (cas formel) 
et 1.2.8 (cas convergent), il y a une 
factorisation $P = Q R$, o\`u $Q$ et $R$ sont entiers unitaires
et $S(Q) = \{\mu\} \;,\; S(R) = S(P) - \{\mu\}$. Le degr\'e de
$Q$ est $r_{P}(\mu) = r_{M}(\mu)$. Le sous-module $M_{Q}$ de
$M = M_{P}$ (2.2.2) est donc pur de pente $\mu$ et de rang
$r_{M}(\mu)$. C'est donc $M^{[\mu]}$ et ce dernier a bien le
rang maximum dans ce cas. \\

On prend maintenant $M$ quelconque. Il existe un entier naturel
non nul $l$ tel que $S(M) \subset \frac{1}{l} \mathbf{Z}$. 
Notons $(K',\sigma_{q'})$ l'extension $(K_{l},\sigma_{q_{l}})$ 
d\'ej\`a utilis\'ee en 2.2.3. Soit $M' = K' \otimes_{K} M$
le module obtenu par extension des scalaires et notons
$\mu' = l \mu$: c'est une pente de $M'$, la plus grande dans 
le cas convergent; de plus, $S(M') = l S(M) \subset \mathbf{Z}$.
D'apr\`es le premier cas \'etudi\'e, ${M'}^{[\mu']}$ a le rang maximum, 
soit $r_{M'}(\mu') = r_{M}(\mu)$. 
Le groupe de Galois de $K'$ sur $K$ est cyclique, engendr\'e
par l'automorphisme $\gamma: f(z') \mapsto f(jz')$, o\`u l'on
a not\'e $z'$ la variable ramifi\'ee $z_{l}$ et o\`u $j$ est
une racine primitive $l$-\`eme de l'unit\'e. L'automorphisme
$K$-lin\'eaire $\gamma \otimes Id_{M}$ de $M'$ commute \`a
$\sigma'$, c'est donc un automorphisme de module aux 
$q'$-diff\'erences,qui laisse ${M'}^{[\mu]}$ invariant
d'apr\`es les remarques pr\'ec\'edentes. Par descente
galoisienne (cf. par exemple \cite{Springer}, chap. 11.1),
il existe un unique sous-espace vectoriel $M_{1}$ du $K$-espace 
vectoriel $M$ tel que ${M'}^{[\mu']} = K' \otimes_{K} M_{1}$.
Il est alors facile de voir que $M_{1}$ est en fait le 
sous-module $M^{[\mu]}$ et que celui-ci a bien pour rang
le rang maximum $r_{M}(\mu) = r_{M'}(\mu')$. 
\hfill $\Box$ \\

De l'additivit\'e du polygone de Newton pour les suites exactes
on tire alors le \\

\textsl{3.1.2 Corollaire. -}
\emph{Le sous-module $M^{[\mu]}$ est tel que
$S(M/M^{[\mu]}) = S(M) - \{\mu\}$.}
\hfill $\Box$ \\

Le lemme d'Adams permet donc de \emph{casser les pentes
dans le cas convergent}. \\

\textsl{3.1.3 Remarque: $q^{-1}$-diff\'erences. -}
Ici, $|q| > 1$ (cf. la remarque 1.2.9) et l'ordre des pentes intervient.
Le module $(M,\Phi)$ sur le corps aux diff\'erences $(K,\sigma_{q})$ 
permet de d\'efinir le module $(M,\Phi^{-1})$ sur le corps aux diff\'erences
$(K,\sigma_{p})$, o\`u $p = q^{-1}$. Celui-ci mod\'elise l'\'equation (1),
vue comme \'equation aux $p$-diff\'erences. Les pentes de ces deux modules
(et de ces deux \'equations) sont deux \`a deux oppos\'ees. Avec les
notations du th\'eor\`eme 3.1.1, le module aux $p$-diff\'erences 
$(M,\Phi^{-1})$ admet un sous-module pur associ\'e \`a la pente $-\mu$
(donc, dans le cas convergent, sa \emph{plus petite} pente); ce sous-module
est form\'e des m\^emes \'el\'ements que $M^{[\mu]}$. \\

\textsl{3.1.4 Un exemple scind\'e. -}
Nous consid\'erons l'op\'erateur:
$$
P = (z \sigma - 1)(\sigma - 1) = z \sigma^{2} - (1+z) \sigma + 1.
$$
Les pentes sont $0$ et $-1$. La factorisation ci-dessus fournit un 
sous-module de pente $0$ de $M = M_{P}$ et le th\'eor\`eme 3.1.1 un 
sous-module de pente $-1$. On a donc $M = M^{[0]} \oplus M^{[-1]}$. 
On le v\'erifie en fabriquant un syst\`eme fondamental de solutions.
On r\'esoud tout d'abord $(\sigma - 1)f = 0$, qui donne $f_{1} = 1$;
puis le syst\`eme $(\sigma - 1)f = g \;,\; (z \sigma - 1)g = 0$,
qui donne $g(z) = e_{q,z^{-1}}$ et 
$f_{2} = \underset{k \geq 1}{\sum} \sigma_{q}^{-k} g$, qui converge. \\

\textsl{3.1.5 Un exemple avec solution divergente. -}
Nous consid\'erons l'op\'erateur:
$$
P = (\sigma - 1)(z \sigma - 1) = q z \sigma^{2} - (1+z) \sigma + 1.
$$
Les pentes sont $0$ et $-1$. La factorisation ci-dessus, qui est la
factorisation canonique, fournit un sous-module de pente $-1$ de 
$M = M_{P}$, en accord avec le th\'eor\`eme 3.1.1: c'est $M^{[-1]}$;
le quotient est pur de pente $0$ (c'est l'unit\'e $\underline{1}$),
mais la suite exacte n'est pas scind\'ee. On le v\'erifie en fabriquant 
un syst\`eme fondamental de solutions. On r\'esoud tout d'abord 
$(z \sigma - 1)f = 0$, qui donne $f_{1}(z) = e_{q,z^{-1}}$.
puis le syst\`eme $(z \sigma - 1)f = g \;,\; (\sigma - 1)g = 0$,
qui donne $g(z) = 1$ et 
$f_{2}(z) = - \underset{k \geq 0}{\sum} q^{k(k-1)/2} z^{k}$, qui diverge
(c'est un $q$-analogue de la s\'erie d'Euler).

\subsubsection*{La filtration canonique}

Soit $M$ un module aux $q$-diff\'erences. On num\'erote ses
pentes: $S(M) = \{\mu_{1},\ldots,\mu_{k}\}$. Dans le cas formel,
l'ordre est arbitraire; dans le cas convergent, on suppose que
$\mu_{1} > \cdots > \mu_{k}$. \\

\textsl{3.1.6 Th\'eor\`eme. -}
\emph{Dans le cas convergent, Il existe une unique tour de
sous-modules:
$\{0\} = M_{0} \subset M_{1} \subset \cdots \subset M_{k} = M$
telle que, pour $1 \leq i \leq k$, le module quotient
$M_{i}/M_{i-1}$ est pur de pente $\mu_{i}$. Les rangs de ces
quotients sont alors les $r_{M}(\mu_{i})$.} \\

\Pr
Cela vient tout seul en it\'erant 3.1.1 et 3.1.2.
\hfill $\Box$ \\

\textsl{3.1.7 Th\'eor\`eme. -}
\emph{Dans le cas formel, $M$ admet une unique d\'ecomposition
en somme directe de modules purs. Ceux-ci sont purs de pentes
$\mu_{1},\ldots,\mu_{k}$ et leurs rangs sont les 
$r_{M}(\mu_{i})$.} \\

\Pr
Cela d\'ecoule de 3.1.1 et de 2.3.5.
\hfill $\Box$ \\

\textsl{3.1.8 La filtration canonique: notations. -}
Tout ce qui suit est cons\'equence triviale des deux r\'esultats
pr\'ec\'edents et vise seulement \`a les mettre en forme pour la
suite. Soit $\mu \in \mathbf{Q}$ quelconque (pas n\'ecessairement 
une pente de $M$). \\

On note $M^{\geq \mu}$ ou $F^{\geq\mu} (M)$ 
(resp. $M^{\gt \mu}$ ou $F^{\gt\mu} (M)$)
le plus grand sous-module de $M$ dont toutes les pentes
sont $\geq \mu$ (resp. $\gt \mu$). On a donc:
$$
S(M^{\geq \mu}) = S(M) \; \cap \; [\mu ; +\infty[
\text{ et son rang vaut } \sum_{\mu_{i} \geq \mu} r_{M}(\mu_{i}),
$$
$$
S(M^{\gt \mu}) = S(M) \; \cap \; ]\mu ; +\infty[
\text{ et son rang vaut } \sum_{\mu_{i} \gt \mu} r_{M}(\mu_{i}).
$$
Les $F^{\geq \mu} (M)$ forment la \emph{filtration canonique
(descendante) par les pentes de $M$}, dont les propri\'et\'es
seront \'etudi\'ees en 3.2. \\

On note $M^{(\mu)} = M^{\geq \mu}/M^{\gt \mu}$: c'est $\{0\}$
si $\mu$ n'est pas une pente de $M$, un module pur de pente
$\mu$ et de rang $r_{M}(\mu)$ si $\mu$ est une pente de $M$. \\

On note enfin $Q^{\lt \mu} (M)$ (resp. $Q^{\leq \mu} (M)$ le quotient 
$M/F^{\geq \mu} (M)$ (resp. $M/F^{\gt \mu} (M)$).  
Les couples $(F^{\geq\mu} (M),Q^{\lt \mu} (M))$  et
$(F^{\gt\mu} (M),Q^{\leq \mu} (M))$ cassent donc en deux 
le polygone de Newton.

\subsubsection*{Pentes et croissance des it\'er\'es}

Soit $(M,\Phi)$ un module aux $q$-diff\'erences. 
Il est possible de caract\'eriser les \'el\'ements $x \in M^{\geq \mu}$ 
en termes de vitesse de croissance $z$-adique de la suite des it\'er\'es 
de $x$ par $\Phi$, dans l'esprit du crit\`ere de Jurkat
(\emph{voir} \cite{Katz1}, 11.6). \\

\textsl{3.1.9 Valuations et r\'eseaux. -}
Rappelons que $\mathcal{O}$ d\'esigne l'anneau de valuation 
du corps valu\'e $K$ (pour la valuation $v_{0}$) et que son corps
r\'esiduel est $\frac{\mathcal{O}}{z \mathcal{O}} = \mathbf{C}$. 
Soit $\Lambda$ un r\'eseau du $K$-espace vectoriel de dimension finie $V$
(autrement dit, $\Lambda$ est un sous $\mathcal{O}$-module libre de $V$
de rang $\dim_{K}(V)$). Pour tout $x \in V$, notons:
\begin{eqnarray*}
v_{\Lambda}(x) & = & \sup \{k \in \mathbf{Z} \;/\; x \in z^{k} \Lambda\} \\
               & = & \min (v_{0}(x_{1}),\ldots,v_{0}(x_{n})),
\end{eqnarray*}
o\`u $x_{1},\ldots,x_{n}$ sont les coordonn\'ees de $x$ dans une base 
quelconque du $\mathcal{O}$-module $\Lambda$. 
Pour tout r\'eseau $\Lambda'$, 
$v_{\Lambda'} - v_{\Lambda}$ est une application born\'ee de 
$V - \{0\}$ dans $\mathbf{Z}$. \\

\textsl{3.1.10 Proposition. -}
\emph{Soit $(M,\Phi)$ un module aux $q$-diff\'erences de premi\`ere pente
$\mu$ et de derni\`ere pente $\nu$, et soit $x$ un vecteur cyclique de $M$. 
Alors, pour tout r\'eseau $\Lambda$ de $M$:
\begin{eqnarray*}
v_{\Lambda}(\Phi^{k}(x)) & = & \mu k + O(1) \;,\; k \in \mathbf{N} \\
v_{\Lambda}(\Phi^{-k}(x)) & = & - \nu k + O(1) \;,\; k \in \mathbf{N}.
\end{eqnarray*}}
\Pr
La remarque 3.1.3, permet de d\'eduire imm\'ediatement la deuxi\`eme 
formule de la premi\`ere. Prouvons celle-ci. 
Quitte \`a ramifier, on peut supposer $\mu$ enti\`ere. Notons
$\Psi = z^{- \mu} \Phi$, de sorte que la premi\`ere pente du
module aux $q$-diff\'erences $(M,\Psi)$ vaut $0$ et que $x$
en est un vecteur cyclique. De la formule g\'en\'erale
$$
\forall \alpha \in K \;,\; \left(\alpha \Phi\right)^{k} =
\alpha \sigma_{q}(\alpha) \cdots \sigma_{q}^{k-1}(\alpha) \Phi^{k},
$$
on tire que 
$v_{\Lambda}(\Psi^{k}(x)) = v_{\Lambda}(\Phi^{k}(x)) - \mu k$.
On est donc ramen\'e \`a montrer que $v_{\Lambda}(\Psi^{k}(x))$
est born\'e lorsque $k$ parcourt $\mathbf{N}$. D'apr\`es 3.1.9, 
cette assertion ne d\'epend pas du choix du r\'eseau $\Lambda$. 
Nous prendrons:
$$
\Lambda = \sum_{i = 0}^{n-1} \mathcal{O} \Psi^{i}(x).
$$
Soit $P = \sigma^{n} + a_{1} \sigma^{n-1} + \cdots + a_{n}$
le polyn\^ome minimal de $x$ (pour $\Psi$). On a, d'apr\`es 2.2.7,
$\underset{1 \leq i \leq n}{\min} v_{0}(a_{i}) = 0$. 
En particulier, $a_{1},\ldots,a_{n} \in \mathcal{O}$ et l'on voit
que $\Psi(\Lambda) \subset \Lambda$. On en d\'eduit que
$\Psi(z \Lambda) \subset z \Lambda$, puis que $\Psi$ induit un
endomorphisme $\overline{\Psi}$ du $\mathbf{C}$-espace vectoriel
$\overline{\Lambda} = \frac{\Lambda}{z \Lambda}$. Le vecteur
$\overline{x} = x \pmod{z \Lambda}$ est cyclique pour $\overline{\Psi}$,
dont le polyn\^ome minimal est
$\overline{P} = \sigma^{n} + a_{1}(0) \sigma^{n-1} + \cdots + a_{n}(0)$.
Par hypoth\`ese, il existe un $a_{i}(0) \not= 0 \;(1 \leq i \leq n)$,
et l'on en d\'eduit facilement que, pour tout entier naturel $k$,
$\overline{\Psi}^{k}(\overline{x}) \not= 0$, autrement dit,
$\Psi^{k}(x) \in \Lambda - z \Lambda$, autrement dit,
$v_{\Lambda}(\Psi^{k}(x)) = 0$. 
\hfill $\Box$ \\

On fixe maintenant $a \in ]0;1[$ (par exemple, inspir\'e de la
th\'eorie des nombres, $a = \frac{1}{|q|}$). A toute valuation
$v_{\Lambda}$ sur un $K$-espace vectoriel $V$ est attach\'ee une 
valeur absolue ultram\'etrique sur $V$ d\'efinie par:
$$
\parallel x\parallel_{\Lambda} = a^{v_{\Lambda}(x)}.
$$
De plus, deux telles valeurs absolues sont \'equivalentes. \\

\textsl{3.1.11 Corollaire. -}
\emph{Soit $(M,\Phi)$ un module aux $q$-diff\'erences et soit $x \in M$.
Alors, pour tout r\'eseau $\Lambda$:
$$
\forall \mu \in \mathbf{Q} \;,\; 
x \in M^{\geq \mu} \Leftrightarrow
\parallel\Phi^{k}(x)\parallel_{\Lambda} = O(a^{\mu k}) 
\;,\; k \in \mathbf{N}.
$$}
\Pr
En effet, $x \in M^{\geq \mu}$ \'equivaut \`a $\mu(x) \geq \mu$.
Or, de la proposition 3.1.10, on tire que 
$\parallel\Phi^{k}(x)\parallel_{\Lambda} = a^{k \mu(x) +O(1)}$.
La conclusion est imm\'ediate.
\hfill $\Box$ \\

On prouve de m\^eme: \\

\textsl{3.1.12 Corollaire. -}
\emph{Soit $(M,\Phi)$ un module aux $q$-diff\'erences et soit $x$ 
un vecteur cyclique de $M$. Alors, le module $M$ est fuchsien si
et seulement si $\Phi^{k}(x)$ est born\'e lorsque $k$ parcourt $\mathbf{Z}$
(au sens de n'importe laquelle des valeurs absolues ci-dessus). 
De plus, dans ce cas, pour tout $y \in M$, $\Phi^{k}(y)$ est born\'e 
lorsque $k$ parcourt $\mathbf{Z}$.}
\hfill $\Box$

% 3.2

\subsection{Propri\'et\'es fonctorielles, 
ab\'eliennes et tensorielles}

Tout ce qui suit repose sur des raisonnements par ``abstract
nonsense'' \`a partir du principe suivant, qui est une 
paraphrase de 3.1.6: \\

\emph{Soit $(G^{\geq \mu}(M))_{\mu \in \mathbf{Q}}$ une filtration 
descendante  d'un module aux $q$-diff\'erences $M$ telle que les sauts 
non triviaux $G^{\geq \mu}(M)/G^{\gt \mu}(M)$ (n\'ecessairement en 
nombre fini puisque $M$ est de rang fini) soient purs de pente 
l'indice $\mu$ du saut. Alors c'est la filtration canonique}.

\subsubsection*{Propri\'et\'es fonctorielles et ab\'eliennes}

\textsl{3.2.1 Lemme. -}
\emph{
(i) Soit $M'$ un sous-module de $M$. Alors:
$$
\forall \mu \in \mathbf{Q} \;,\;
F^{\geq \mu}(M') = F^{\geq \mu}(M) \cap M'.
$$
(ii) Soit $M''$ un module quotient de $M$. Alors:
$$
\forall \mu \in \mathbf{Q} \;,\;
F^{\geq \mu}(M'') = \text{ image de } F^{\geq \mu}(M) 
\text{ dans } M''.
$$}
\Pr
En effet, dans chaque cas, le membre de droite est le terme
g\'en\'eral d'une filtration \`a laquelle on peut appliquer le 
principe \'enonc\'e plus haut.
\hfill $\Box$ \\

\textsl{3.2.2 Proposition. -}
\emph{Tout morphisme de modules aux $q$-diff\'erences est strict 
relativement aux filtrations canoniques. Autrement dit, si $f$ 
est un morphisme de $M$ dans $N$: 
$$
\forall \mu \in \mathbf{Q} \;,\;
f(M^{\geq \mu}) = f(M) \cap N^{\geq \mu}.
$$}
\Pr
Cela d\'ecoule imm\'ediatement du lemme.
\hfill $\Box$ \\

\textsl{3.2.3 Corollaire. -}
\emph{Les $F^{\geq \mu}$ et les $F^{\gt \mu}$ sont des 
endofoncteurs de la cat\'egorie $DiffMod(K,\sigma_{q})$.}
\hfill $\Box$ \\

Les propri\'et\'es sp\'ecifiquement ab\'eliennes (comportement 
vis \`a vis des suites exactes) viennent de la proposition et
apparaitront plus clairement en 3.3, avec l'\'etude du foncteur 
``gradu\'e associ\'e''.

\subsubsection*{Propri\'et\'es tensorielles}

On part d'un produit tensoriel $M = M_{1} \otimes M_{2}$. 
On a donc $S(M) = S(M_{1}) + S(M_{2})$ (cf. 2.3.3). \\

\textsl{3.2.4 Proposition. -}
\emph{La filtration canonique sur le produit tensoriel
$M_{1} \otimes M_{2}$ est donn\'ee par la formule:
$$
F^{\geq \mu}(M_{1} \otimes M_{2}) =
\sum_{\mu_{1} + \mu_{2} \geq \mu} 
F_{1}^{\geq \mu_{1}}(M_{1}) \otimes F_{2}^{\geq \mu_{2}}(M_{2}).
$$}
\Pr
Pour $\mu \in \mathbf{Q}$, on pose:
$$
G^{\geq \mu}(M) = 
\sum_{\mu_{1} + \mu_{2} \geq \mu} 
M_{1}^{\geq \mu_{1}} \otimes M_{2}^{\geq \mu_{2}}.
$$
Il est clair que les $(G^{\geq \mu}(M))_{\mu \in \mathbf{Q}}$ 
forment une filtration descendante de $M$ dont les sauts ont lieu
aux pentes de $M$. \\

Notons d'autre part:
$$
H^{\mu}(M) = 
\sum_{\mu_{1} + \mu_{2} \geq \mu} 
\left(M_{1}^{\geq \mu_{1}} \otimes M_{2}^{\gt \mu_{2}} +
M_{1}^{\gt \mu_{1}} \otimes M_{2}^{\geq \mu_{2}}\right)
$$
On voit que $H^{\mu}(M) \subset G^{\geq \mu}(M)$, l'inclusion 
\'etant stricte si et seulement si $\mu$ est une pente de $M$.
Soient $\mu' > \mu$ deux pentes cons\'ecutives de $M$.
Alors, si $\mu_{1} + \mu_{2} \geq \mu$, on voit que
$M_{1}^{\gt \mu_{1}} \otimes M_{2}^{\geq \mu_{2}}
\subset H^{\mu'}(M)$ puis, par sym\'etrie et en additionnant,
que $H^{\mu}(M) \subset G^{\geq \mu'}(M)$. Par cons\'equent,
$H^{\mu}(M) = G^{\gt \mu}(M)$. \\

On a deux surjections naturelles:
$$
\bigoplus_{\mu_{1} + \mu_{2} \geq \mu} 
M_{1}^{\geq \mu_{1}} \otimes M_{2}^{\geq \mu_{2}} 
\rightarrow G^{\geq \mu}(M) \rightarrow 0
$$
et
$$
\bigoplus_{\mu_{1} + \mu_{2} \geq \mu} 
\left(M_{1}^{\geq \mu_{1}} \otimes M_{2}^{\gt \mu_{2}} +
M_{1}^{\gt \mu_{1}} \otimes M_{2}^{\geq \mu_{2}}\right)
\rightarrow G^{\gt \mu}(M) \rightarrow 0,
$$
d'o\`u, par passage au quotient, une surjection:
$$
\bigoplus_{\mu_{1} + \mu_{2} \geq \mu} 
\frac{M_{1}^{\geq \mu_{1}} \otimes M_{2}^{\geq \mu_{2}}}
{M_{1}^{\geq \mu_{1}} \otimes M_{2}^{\gt \mu_{2}} +
M_{1}^{\gt \mu_{1}} \otimes M_{2}^{\geq \mu_{2}}}
\rightarrow \frac{G^{\geq \mu}(M)}{G^{\gt \mu}(M)}
\rightarrow 0.
$$
Or, de la formule g\'en\'erale :
$$
\frac{A'}{A} \otimes \frac{B'}{B} =
\frac{A' \otimes B'}{A \otimes B' + A' \otimes B},
$$
on d\'eduit l'isomorphisme canonique:
$$
\frac{M_{1}^{\geq \mu_{1}} \otimes M_{2}^{\geq \mu_{2}}}
{M_{1}^{\geq \mu_{1}} \otimes M_{2}^{\gt \mu_{2}} +
M_{1}^{\gt \mu_{1}} \otimes M_{2}^{\geq \mu_{2}}} =
M_{1}^{(\mu_{1})} \otimes M_{2}^{(\mu_{2})},
$$
autrement dit, on a une surjection:
$$
\bigoplus_{\mu_{1} + \mu_{2} = \mu} 
M_{1}^{(\mu_{1})} \otimes M_{2}^{(\mu_{2})}
\rightarrow \frac{G^{\geq \mu}(M)}{G^{\gt \mu}(M)}
\rightarrow 0,
$$
qui entraine que $G^{\geq \mu}(M)/G^{\gt \mu}(M)$
est pur de pente $\mu$.
\hfill $\Box$ \\

\textsl{3.2.5 Remarque. -}
Ces propri\'et\'es de la filtration sont analogues \`a celles
axiomatis\'ees par Saavedra dans dans \cite{Saavedra} (en mieux 
puisque nous filtrons par des sous-objets) et \emph{oppos\'ees} 
\`a celles axiomatis\'ees par Yves Andr\'e dans \cite{YA3}.

% 3.3

\subsection{Le gradu\'e associ\'e}

On r\'ecolte ici presque gratuitement les fruits du travail
d\'ej\`a fait.

\subsubsection*{Modules mod\'er\'ement irr\'eguliers}

\textsl{3.3.1 La cat\'egorie $DiffMod_{mi}(K,\sigma_{q})$. -} 
On appellera \emph{mod\'er\'ement irr\'egulier} un module
aux $q$-diff\'erences somme directe de modules purs. Cette
terminologie est inspir\'ee de l'expression ``mod\'er\'ement
ramifi\'e'' en arithm\'etique et en g\'eom\'etrie alg\'ebrique.
La propri\'et\'e d'\^etre mod\'er\'ement irr\'egulier est 
conserv\'ee par extension des scalaires. 
On notera $DiffMod_{mi}(K,\sigma_{q})$ la sous-cat\'egorie 
pleine de $DiffMod(K,\sigma_{q})$ form\'ee des modules 
mod\'er\'ement irr\'eguliers et, pour tout entier naturel non nul 
$l$ (niveau de ramification suffisant pour triangulariser), 
on notera $DiffMod_{mi,l}(K,\sigma_{q})$ la sous-cat\'egorie
pleine dont les objets ont toutes leurs pentes dans 
$\frac{1}{l} \mathbf{Z}$. 
La cat\'egorie $DiffMod_{mi}(K,\sigma_{q})$ est naturellement
gradu\'ee par $\mathbf{Q}$. Les morphismes sont de degr\'e $0$.
La cat\'egorie est stable par sous-quotient et constructions
tensorielles (mais pas par extensions). C'est donc une 
sous-cat\'egorie tannakienne de $DiffMod(K,\sigma_{q})$.
Ces propri\'et\'es sont encore valables pour chacune des
sous-cat\'egories $DiffMod_{mi,l}(K,\sigma_{q})$. 
La projection $M \leadsto M^{(\mu)}$ sur la composante de
degr\'e (ou pente) $\mu$ est un foncteur exact. \\

\textsl{3.3.2 Fonction et polyn\^ome de Hilbert-Samuel. -}
On peut associer \`a tout module mod\'er\'ement irr\'egulier
$M = \bigoplus M^{(\mu)}$ sa \emph{fonction de Hilbert-Samuel}
$\mu \mapsto r_{M}(\mu) = \dim_{K}(M^{(\mu)})$ (qui est aussi
sa fonction de Newton) et son \emph{polyn\^ome de Hilbert-Samuel}:
$$
\mathcal{R}_{M}(T) = 
\sum_{\mu \in \mathbf{Q}} r_{M}(\mu) T^{\mu},
$$
(qui est un polyn\^ome ramifi\'e). La formation de ce dernier 
est additive pour les suites exactes et multiplicative pour le 
produit tensoriel. La dualisation se traduit par 
$T \leftarrow T^{-1}$. En combinaison avec 3.3.4,
ces propri\'et\'es ``expliquent'' (ou enrichissent) les 
propri\'et\'es de la fonction de Newton obtenues en 2.3.1.

\subsubsection*{Propri\'et\'es fonctorielles,ab\'eliennes
et tensorielles}

Si $M$ est un module aux $q$-diff\'erences, rappelons que 
l'on note $M^{(\mu)} = \frac{F^{\geq \mu}(M)}{F^{\gt \mu}(M)}$
le facteur de pente $\mu$ de $M$ et que celui-ci est non nul
si et seulement si $\mu$ est une pente de $M$; dans ce cas,
il est pur de pente $\mu$ et de rang $r_{M}(\mu)$. Nous noterons: 
$$
gr(M) = \bigoplus_{\mu \in \mathbf{Q}} M^{(\mu)}
      = \bigoplus_{\mu \in S(M)} M^{(\mu)}
$$
le \emph{gradu\'e associ\'e} \`a la filtration canonique
de $M$; c'est un module mod\'er\'ement irr\'egulier. 
D'apr\`es 3.2.2, tout morphisme $f: M \rightarrow N$ envoie
$F^{\geq \mu}(M)$ dans $F^{\geq \mu}(N)$ et $F^{\gt \mu}(M)$
dans $F^{\gt \mu}(N)$, donc induit 
$f^{(\mu)}: M^{(\mu)} \rightarrow N^{(\mu)}$
donc aussi $gr(f): gr(M) \rightarrow gr(N)$. \\

\textsl{3.3.3 Th\'eor\`eme. -}
\emph{On d\'efinit ainsi un foncteur $\mathbf{C}$-lin\'eaire
exact de la cat\'egorie ab\'elienne $DiffMod(K,\sigma_{q})$
dans sa sous-cat\'egorie pleine $DiffMod_{mi}(K,\sigma_{q})$.
Ce foncteur est une r\'etraction de l'inclusion.} \\

\Pr
L'exactitude est une cons\'equence classique du fait que tous
les morphismes sont stricts (cf 3.2.2).
\hfill $\Box$ \\

\textsl{3.3.4 Th\'eor\`eme. -}
\emph{Le foncteur $gr$ est compatible au produit tensoriel
et fid\`ele.} \\

\Pr
D'apr\`es \cite{Saavedra}, IV.2.1, il d\'ecoule de 3.2.4
que le foncteur qui associe \`a $M$ le $K$-espace vectoriel
sous-jacent \`a $gr(M)$ poss\`ede les propri\'et\'es
indiqu\'ees. Le fait que $gr$ lui-m\^eme poss\`ede ces
propri\'et\'es est alors trivial.
\hfill $\Box$ \\

\textsl{3.3.5 Remarque. -}
D'apr\`es 3.1.7, dans le cas formel, 
$DiffMod_{mi}(K,\sigma_{q}) = DiffMod(K,\sigma_{q})$
et ce foncteur est isomorphe au foncteur identit\'e.

% 3.4

\subsection{Applications \`a la classification et 
\`a la th\'eorie de Galois}

Il s'agit ici d'esquisses; pour les preuves d\'etaill\'ees, 
voir les r\'ef\'erences indiqu\'ees.

\subsubsection*{Applications \`a la classification}

La filtration par les pentes et le gradu\'e associ\'e sont 
un des outils de la \emph{classification analytique locale 
par voie transcendante des \'equations aux $q$-diff\'erences
lin\'eaires \`a coefficients rationnels}. Celle-ci est expos\'ee 
dans l'article en pr\'eparation \cite{RSZ}. La partie alg\'ebrique,
que nous r\'esumons ici (c'est la plus facile), est accessible dans 
\cite{JSGTQDIFCAL}. \\

D'apr\`es 3.3, la classification formelle se ram\`ene \`a celle des
modules purs et la classification locale ``analytique modulo 
formelle'' se ram\`ene \`a la classification \`a gradu\'e donn\'e.
Soient $P_{1},\ldots,P_{k}$ des modules purs de rangs $r_{1},\ldots,r_{k}$ 
et de pentes $\mu_{1} \gt \cdots \gt \mu_{k}$. On pose:
$$
\mathcal{F}(P_{1},\ldots,P_{k}) =
\{\text{classes d'isomorphies de couples } \; (M,g)\},
$$
o\`u les couples $(M,g)$ sont form\'es d'un module $M$ et d'un isomorphisme 
$g : gr(M) \rightarrow P_{1} \oplus \cdots \oplus P_{k}$ et $(M,g)$ est 
isomorphe \`a $(M',g')$ s'il existe un isomorphisme de modules 
$u : M \rightarrow M'$ tel que $g = g' \circ gr(u)$. La description de 
$\mathcal{F}(P_{1},\ldots,P_{k})$ se ram\`ene \`a des calculs de classes 
d'extensions. On prouve:

\begin{enumerate}

\item{Que le foncteur ``sections globales'' $\Gamma$ est exact \`a gauche,
que les $\Gamma^{i}(M) = Ext^{i}(\underline{1},M)$ sont les foncteurs 
d\'eriv\'es du foncteur $\Gamma$ et que 
$Ext^{i}(M,N) = \Gamma^{i}(M^{\vee} \otimes N)$.}

\item{Que les $Ext^{i} , i \geq 2$ sont nuls et que les 
$Ext^{i} , i = 0 , 1$ sont de dimension finie.}

\end{enumerate}

On d\'efinit alors la \emph{caract\'eristique d'Euler-Poincar\'e}
$\chi(M) = \dim \Gamma^{0}(M) - \dim \Gamma^{1}(M)$, qui est
additive pour les suites exactes. Dans le cas d'objets de rang
$1$, $\Gamma^{0}$ et $\Gamma^{1}$ s'interpr\`etent
respectivement comme un noyau et un conoyau d'op\'erateur aux 
$q$-diff\'erences et la caract\'eristique d'Euler-Poincar\'e
comme un \emph{indice}, ce qui permet des calculs exacts. 
On peut alors, en analogie avec \cite{Birkhoff3}, donner \`a
un module \`a pentes enti\`eres une forme normale polynomiale
et en d\'eduire, pour $\mathcal{F}(P_{1},\ldots,P_{k})$,
une structure de vari\'et\'e alg\'ebrique affine de dimension
$\underset{1 \leq i < j \leq k}{\prod} r_{i} r_{j} (\mu_{i} - \mu_{j})$. \\

Une nouvelle m\'ethode de resommation discr\`ete permet d'uniformiser 
cette vari\'et\'e \`a l'aide de fonctions elliptiques \`a p\^oles contr\^ol\'es
et une th\'eorie adapt\'ee des d\'eveloppements asymptotiques en fournit
une interpr\'etation en termes du faisceau de Malgrange, ici d\'efini
sur la courbe elliptique $\mathbf{E}_{q}$.

\subsubsection*{Applications \`a la th\'eorie de Galois}

Ces r\'esultats figurent dans l'article en pr\'eparation \cite{JSIRR}.
Ils concernent le cas convergent, le cas formel \'etant plus simple
(\emph{voir} \cite{SVdP}). On se restreint ici au cas ``non ramifi\'e''.
Notons $DiffMod_{1}(K,\sigma_{q})$ la sous-cat\'egorie
pleine de $DiffMod(K,\sigma_{q})$ form\'ee des objets \`a pentes 
enti\`eres. C'est une sous-cat\'egorie tannakienne, stable
par extensions. On d\'eduit ais\'ement de \cite{JSGAL} que le groupe de
Galois $G_{mi,1}^{(0)}$ de $DiffMod_{mi,1}(K,\sigma_{q})$
est \'egal \`a $\mathbf{C}^{*} \times G_{f}^{(0)}$, o\`u
$G_{f}^{(0)}$ est le groupe de Galois local fuchsien, 
qui est isomorphe \`a:
$$
Hom_{gr}(\mathbf{C}^{*}/q^{\mathbf{Z}},\mathbf{C}^{*})
\times \mathbf{C}.
$$
D'apr\`es 3.3.3,3.3.4 et \cite{DM}, le groupe de Galois $G_{irr,1}^{(0)}$
de $DiffMod_{1}(K,\sigma_{q})$ est le produit semi-direct 
d'un groupe unipotent et de $G_{mi,1}^{(0)}$.
Le groupe unipotent est form\'e des automorphismes galoisiens
dont l'effet ne se voit pas sur le gradu\'e, c'est \`a dire
ceux qui ont vocation \`a \^etre consid\'er\'es comme des Stokes.  
Conjecturalement, il y a une description en termes de fibr\'es
sur la courbe elliptique pour laquelle la filtration canonique
correspond \`a la filtration de Harder-Narasimhan.

\appendix

% A

\section{Solutions formelles et solutions convergentes}

Dans toute cette section, les pentes $\mu$ seront des entiers.
D'apr\`es ce qui pr\'ec\`ede, nous sommes conduits
\`a nous int\'eresser \`a l'\'equation avec second
membre: 
$$
z^{-\mu} \sigma_{q} f - c f = g
$$ 
laquelle se ram\`ene, apr\`es transformation de jauge de symbole 
$e_{q,c z^{\mu}}$ (1.1.7), \`a l'\'equation:
$$
\sigma_{q} f - f = g.
$$ 
Si l'on adopte l'analogie habituelle avec le cas 
diff\'erentiel:
$$
\frac{\sigma_{q} - 1}{q - 1} 
\longleftrightarrow z \frac{d}{dz},
$$
on est conduit \`a consid\'erer cette r\'esolution
comme une \emph{$q$-int\'egration}. Comme dans le
cas diff\'erentiel, la constante $1$ n'est pas $q$-int\'egrable 
et n\'ecessite l'introduction du $q$-logarithme.

\subsection*{A.1~~~$q$-int\'egration}

Soit $\pi_{0}$ le projecteur du $\mathbf{C}$ espace vectoriel 
$\mathbf{C}((z))$ qui associe \`a toute s\'erie de Laurent 
formelle son terme constant. Les sous-espaces vectoriels
$\mathcal{M}(\mathbf{C})$ et $\mathbf{C}(\{z\})$ sont stables,
d'o\`u, quelque soit le corps $K$, une d\'ecomposition:
$$
K = \mathbf{C} \oplus K^{\bullet},
\text{ o\`u }
K^{\bullet} = \left(\Ker \pi_{0}\right) \cap K.
$$
L'endomorphisme $\mathbf{C}$-lin\'eaire $\sigma_{q} - 1$ 
de $K$ est nul sur la premi\`ere composante et laisse stable
la seconde. \\

\textsl{A.1.1 Lemme. -}
\emph{L'endomorphisme $\sigma_{q} - 1$ induit un automorphisme 
de $K^{\bullet}$.} \\

\Pr
En effet, on peut poser (dans $\mathbf{C}((z))^{\bullet}$):
$$
I_{q} \left(\sum_{i \not= 0} a_{i} z^{i}\right) =
\sum_{i \not= 0} \frac{a_{i}}{q^{i}-1} z^{i},
$$
d\'efinissant un inverse. Il est clair que celui-ci pr\'eserve,
le cas \'ech\'eant, la m\'eromorphie pr\`es de $0$ ou
sur $\mathbf{C}$.
\hfill $\Box$ \\

On introduit donc maintenant un \'el\'ement $l_{q}$ de $L$
tel que $\sigma_{q} l_{q} = l_{q} + 1$ (voir dans l'introduction 
les conventions g\'en\'erales). On note de plus, pour tout
entier naturel $k$:
$$
l_{q}^{(k)} = \begin{pmatrix} l_{q} \\ k \end{pmatrix} =
\frac{1}{k!} \prod_{i=0}^{k-1} (l_{q} - i),
$$
et $l_{q}^{(k)} = 0$ pour $k \lt 0$, de sorte que 
(calcul facile):
$$
\forall k \in \mathbf{Z} \;,\;
\sigma_{q} l_{q}^{(k)} = l_{q}^{(k)} + l_{q}^{(k-1)}.
$$

\textsl{A.1.2 Lemme. -}
\emph{Les $l_{q}^{(k)}$, $k \geq 0$, sont lin\'eairement 
ind\'ependants sur $K$; autrement dit, $l_{q}$ est
transcendant et :
$$
K[l_{q}] = \bigoplus_{k \geq 0} K l_{q}^{(k)}.
$$}
\Pr
Soit en effet une relation:
$$
l_{q}^{(k+1)} = 
a_{0} \; l_{q}^{(0)} + \cdots + a_{k} \; l_{q}^{(k)},
\quad \text{les } a_{i} \in K,
$$
avec $k \geq 0$ le plus petit possible; il est donc en fait
$\geq 1$ puisque $l_{q} \not\in K$. En appliquant $\sigma_{q} - 1$, 
\`a cette relation, on trouve:
$$
l_{q}^{(k)} \equiv (\sigma_{q} a_{k} - a_{k}) l_{q}^{(k)}
\pmod{K l_{q}^{(0)} + \cdots + K l_{q}^{(k-1)}}.
$$
Par minimalit\'e, on en d\'eduit que
$\sigma_{q} a_{k} - a_{k} = 1$, ce qui est impossible.
\hfill $\Box$ \\

\textsl{A.1.3 Proposition. -}
\emph{On a, pour tout entier naturel non nul $k$, 
une suite exacte:
$$
0 \rightarrow \mathbf{C} \rightarrow K_{k}[l_{q}]
\overset{\sigma_{q} - 1}{\longrightarrow} K_{k-1}[l_{q}] 
\rightarrow 0.
$$}
\Pr
Ici, $K_{k}[X]$ d\'esigne l'ensemble des polyn\^omes de
degr\'e $\leq k$.
Ecrivons $f = f_{0} l_{q}^{(0)} + \cdots + f_{k} l_{q}^{(k)}$
et $g = g_{0} l_{q}^{(0)} + \cdots + g_{k-1} l_{q}^{(k-1)}$
des \'el\'ements respectifs de $K_{k}[l_{q}]$ et de
$K_{k-1}[l_{q}]$. Par identification, l'\'equation
$(\sigma_{q} - 1) f = g$ \'equivaut \`a:
$$
\forall i \geq 0 \;,\; 
g_{i} = \sigma_{q} f_{i} - f_{i} + \sigma_{q} f_{i+1}.
$$
La r\'esoudre revient \`a r\'esoudre le syst\`eme:
$$
\begin{cases}
\sigma_{q} f_{0} - f_{0} + \sigma_{q} f_{1} = g_{0} \\
\vdots \\
\sigma_{q} f_{i} - f_{i} + \sigma_{q} f_{i+1} = g_{i} \\
\vdots \\
\sigma_{q} f_{k-1} - f_{k-1} + \sigma_{q} f_{k} = g_{k-1} \\
\sigma_{q} f_{k} - f_{k} = 0
\end{cases}
$$
On voit, en commen\c{c}ant par le bas, que $f_{k} \in \mathbf{C}$
et m\^eme (avant-derni\`ere \'equation) que c'est n\'ecessairement
$\pi_{0}(g_{k-1})$. On a alors la r\'esolution it\'erative:
\begin{eqnarray*}
f_{k} & = & \pi_{0}(g_{k-1}) \\
      & \vdots &             \\
f_{i} & = & \pi_{0}(g_{i-1}) + I_{q}(g_{i} - \sigma_{q} f_{i+1}) \\
      & \vdots &             \\
f_{0} & = & \text{ une constante arbitraire } + 
            I_{q}(g_{0} - \sigma_{q} f_{1})
\end{eqnarray*}
\hfill $\Box$ 

\subsection*{A.2~~~Equations d'ordre $1$ avec second membre}

On se restreint dor\'enavant \`a la sous-alg\`ebre $S$ de $L$
engendr\'ee par les fonctions \'el\'ementaires:
$$
S = K [(e_{q,cz^{\mu}})_{(c,\mu) \in \mathbf{C}^{*} \times \mathbf{Z}},l_{q}].
$$
Notons provisoirement $C(S)$ l'ensemble des ``caract\`eres'':
$$
C(S) = \{u \in S - \{0\} \;/\; 
         \exists c \in \mathbf{C}^{*} \;:\; \sigma_{q} u = c u\}.
$$
On a une d\'ecomposition
\footnote{Cette d\'ecomposition poss\`ede d'int\'eressantes 
propri\'et\'es alg\'ebriques, partiellement abord\'ees dans
\cite{SVdP} (cas formel) et \cite{JSAIF} (cas convergent).}
:
$$
S = \sum_{\mu \in \mathbf{Z}} S_{\mu}, \quad \text{o\`u} \quad
S_{\mu} = \sum_{u \in C(S)} u \Theta_{q}^{\mu} K[l_{q}].
$$

La formule, imm\'ediatement v\'erifi\'ee:
$$
\sigma_{q} u = c u \Rightarrow
(d z^{\nu} \sigma -1)(u \Theta_{q}^{\mu} F) =
u \Theta_{q}^{\mu} (cd z^{\mu + \nu} \sigma -1) F
$$
implique que l'endomorphisme $\Phi_{d,\nu} = d z^{\nu} \sigma -1$ 
du $\mathbf{C}$-espace vectoriel $S$ laisse stable chaque
sous-espace $u \Theta_{q}^{\mu} K[l_{q}]$. De plus, 
l'isomorphisme $F \mapsto u \Theta_{q}^{\mu} F$ de
$K[l_{q}]$ dans $u \Theta_{q}^{\mu} K[l_{q}]$ conjugue
l'action de $\Phi_{cd,\mu+\nu}$ sur le premier avec
l'action de $\Phi_{d,\nu}$ sur le deuxi\`eme.
Notre but, dans ce paragraphe, est de pr\'eciser l'image
et le noyau de ces endomorphismes, et, en particulier,
de d\'emontrer le th\'eor\`eme A.2.4. \\

\textsl{A.2.1 Lemme. -}
\emph{Soit $(c,\mu) \in \mathbf{C}^{*} \times \mathbf{Z}$.
Il est clair que $K$ est stable par $\Phi_{c,\mu}$. \\
(i) Si $(\overline{c},\mu) \not= (1,0)$, la restriction
de $\Phi_{c,\mu}$ \`a $K$ est injective. \\
(ii) Elle est de plus surjective dans chacun des cas
suivants:
\begin{enumerate}
\item{$\mu = 0$ et $\overline{c} \not= 1$.}
\item{$\mu \lt 0$.}
\item{$\mu \gt 0$ et $K = \mathbf{C}((z))$.}
\end{enumerate}
}

\Pr
Si $\mu = 0$, \'ecrivant
$f = \underset{k \gt \gt -\infty}{\sum} f_{k} z^{k}$ et
$g = \underset{k \gt \gt -\infty}{\sum} g_{k} z^{k}$,
on obtient l'\'equivalence:
$$
(c \sigma -1)f = g \Leftrightarrow
\forall k \in \mathbf{Z} \;,\; (c q^{k} -1) f_{k} = g_{k},
$$
qui suffit \`a montrer (i) et (ii) dans ce cas
(c'est l'hypoth\`ese $\overline{c} \not= 1$ 
qui garantit que $c q^{k} - 1$ ne s'annule pas). \\

Si $\mu \not= 0$, posons $\mu = m \epsilon$,
avec $m = |\mu|$ et $\epsilon = \pm 1$.
La d\'ecomposition:
$$
\mathbf{C}((z)) = 
\bigoplus_{0 \leq i \lt m} z^{i} \mathbf{C}((z^{m}))
$$
induit des d\'ecompositions similaires de
$\mathbf{C}(\{z\})$ et de $\mathcal{M}(\mathbf{C})$.
La formule (facile \`a v\'erifier):
$$
(c z^{\mu} \sigma -1) z^{i} F(z^{m}) =
z^{i} (c q^{i} z^{m \epsilon} F(q^{m} z^{m}) - F(z^{m}))
$$
montre que chaque composante est stable. Ecrivant alors
$Z = z^{m}$, $Q = q^{m}$, $C = c q^{i}$,
$f(z) = z^{i} F(z^{m})$ et $g(z) = z^{i} G(z^{m})$, 
on obtient l'\'equivalence:
$$
(c z^{\mu} \sigma -1)f = g \Leftrightarrow
C Z^{\epsilon} F(QZ) - F(Z) = G(Z).
$$
Autrement dit, on s'est ramen\'e au cas o\`u $\mu = \pm 1$,
ce que l'on suppose maintenant. On reprend les notations
$f = \underset{k \gt \gt -\infty}{\sum} f_{k} z^{k}$ et
$g = \underset{k \gt \gt -\infty}{\sum} g_{k} z^{k}$. \\

Si $\mu = -1$, on obtient les \'equivalences:
\begin{eqnarray*}
(c z^{-1}\sigma -1)f = g & \Leftrightarrow &
\forall k \in \mathbf{Z} \;,\; 
c q^{k+1} f_{k+1} - f_{k} = g_{k} \\
                         & \Leftrightarrow &
\forall k \in \mathbf{Z} \;,\; 
c^{k+1} q^{k(k+1)/2} f_{k+1} - 
c^{k} q^{k(k-1)/2} f_{k} = 
c^{k} q^{k(k-1)/2} g_{k} \\
                         & \Leftrightarrow &
\forall k \in \mathbf{Z} \;,\; 
c^{k} q^{k(k-1)/2} f_{k} = 
\sum_{i \lt k} c^{i} q^{i(i-1)/2} g_{i}
\end{eqnarray*}
Ceci montre que $\Phi_{c,-1}$ est bijectif dans 
le cas formel. Dans le cas convergent, la relation
$$
|c^{k} f_{k}| \leq \sum_{i \lt k} |c^{i} g^{i}|
$$
entraine que la s\'erie $f(cz)$ est domin\'ee
par la s\'erie $\frac{g(cz)}{1-z}$, ce qui conclut
encore. \\

Si $\mu = +1$, on obtient les \'equivalences:
\begin{eqnarray*}
(c z \sigma -1)f = g & \Leftrightarrow &
\forall k \in \mathbf{Z} \;,\; 
c q^{k-1} f_{k-1} - f_{k} = g_{k} \\
                         & \Leftrightarrow &
\forall k \in \mathbf{Z} \;,\; 
\frac{f_{k-1}}{(c/q)^{k-1} q^{k(k-1)/2}} - 
\frac{f_{k}}{(c/q)^{k} q^{k(k+1)/2}} = 
\frac{g_{k}}{(c/q)^{k} q^{k(k+1)/2}} \\
                         & \Leftrightarrow &
\forall k \in \mathbf{Z} \;,\; 
\frac{f_{k}}{(c/q)^{k} q^{k(k+1)/2}} = 
- \sum_{i \lt k} \frac{g_{i}}{(c/q)^{i} q^{i(i+1)/2}}
\end{eqnarray*}
Ceci montre que $\Phi_{c,1}$ est bijectif dans 
le cas formel. 
\hfill $\Box$ \\

Dans le cas convergent avec $\mu \gt 0$,
on ne peut pas en g\'en\'eral conclure, 
les coefficients $f_{k}$ pouvant \^etre
tr\`es rapidement croissants. Par exemple, 
si $c = 1$ et $g = -1$, on trouve, pour $k \gt 0$,
$f_{k} = q^{k(k-1)/2}$.
C'est un $q$-analogue de la s\'erie d'Euler. \\

On va maintenant \'etudier l'action de $\Phi_{c,\mu}$
sur $K[l_{q}]$. Le cas o\`u $(c,\mu) = (1,0)$ a fait
l'objet du A.1. Le cas o\`u $(\overline{c},\mu) = (1,0)$
s'y ram\`ene car l'automorphisme $F \mapsto z^{l} F$
de $K[l_{q}]$ conjugue $\Phi_{c,\mu}$ avec
$\Phi_{q^{l} c,\mu}$. \\

\textsl{A.2.2 Corollaire. -}
\emph{On suppose $(\overline{c},\mu) \not= (1,0)$.
Les conclusions sont les m\^emes: la restriction
de $\Phi_{c,\mu}$ \`a $K[l_{q}]$ est injective;
elle est de plus surjective, sauf dans le cas
convergent si $\mu \gt 0$.} \\

\Pr
Ecrivant $f = \underset{i \geq 0}{\sum} f^{(i)} l_{q}^{(i)}$
et $g = \underset{i \geq 0}{\sum} g^{(i)} l_{q}^{(i)}$
(qui sont des sommes finies), on obtient l'\'equivalence:
$$
(c z^{\mu} \sigma -1)f = g \Leftrightarrow
\forall i \geq 0 \;,\; 
(c z^{\mu} \sigma - 1) f^{(i)} =
g^{(i)} - c z^{\mu} \sigma_{q} f^{(i+1)}.
$$
Ce syst\`eme se r\'esoud it\'erativement, en commen\c{c}ant
par la fin, \`a l'aide du lemme A.2.1.
\hfill $\Box$ \\

\textsl{A.2.3 Corollaire. -}
\emph{On consid\`ere la restriction de $\Phi_{d,\nu}$ 
\`a $u \Theta_{q}^{\mu} K[l_{q}]$, o\`u 
$\sigma_{q} u = c u$. \\
(i) Si $cd = q^{l} , l \in \mathbf{Z}$,
et si $\mu + \nu = 0$, cet endomorphisme est surjectif 
de noyau $\mathbf{C} u \Theta_{q}^{\mu} z^{-l}$. \\
(ii) Si $\overline{cd} \not= 1$ et $\mu + \nu = 0$,
ou bien si $cd$ est quelconque et $\mu + \nu \lt 0$,
l'endomorphisme est bijectif. \\
(iii) M\^eme conclusion dans le cas formel si
$\mu + \nu \gt 0$.} \\

\Pr
C'est imm\'ediat par conjugaison (voir le d\'ebut
de A.2). 
\hfill $\Box$ \\

Nous synth\'etisons maintenant les r\'esultats les
plus importants: \\

\textsl{A.2.4 Th\'eor\`eme. -}
\emph{L'endomorphisme $\Phi_{d,\nu}$ de $S_{\mu}$
est surjectif si $\mu + \nu \leq 0$, et aussi
si $\mu + \nu \gt 0$ dans le cas formel.}
\hfill $\Box$

\subsection*{A.3~~~R\'esolution formelle}

\textsl{A.3.1 D\'efinition. -}
\emph{Soient $f_{1},\ldots,f_{m}$ des \'el\'ements de $L$.
Leur $q$-Wronskien (ou Casoratien, ou Pochhammerien) est:
$$
W_{q}(f_{1},\ldots,f_{m}) =
\det \begin{pmatrix}
f_{1} & \ldots & f_{j} & \ldots & f_{m} \\
\vdots & \vdots & \vdots & \vdots & \vdots \\
\sigma_{q}^{i} f_{1} & \ldots & \sigma_{q}^{i} f_{j} & 
\ldots & \sigma_{q}^{i} f_{m} \\
\vdots & \vdots & \vdots & \vdots & \vdots \\
\sigma_{q}^{m-1} f_{1} & \ldots & \sigma_{q}^{m-1} f_{j} & 
\ldots & \sigma_{q}^{m-1} f_{m}
\end{pmatrix}
$$}

Rappelons (cf. l'introduction) que l'on note $C_{L} = L^{\sigma_{q}}$ 
le sous-corps des constantes de $L$. Dans ces conditions, on a le: \\

\textsl{A.3.2 Lemme. -}
\emph{Le $q$-Wronskien $W_{q}(f_{1},\ldots,f_{m})$ est non nul
si et seulement si les $f_{i}$ sont lin\'eairement ind\'ependants
sur $C_{L}$.} \\

\Pr
Ce lemme est d\'emontr\'e dans \cite{LDVpreprint}.
\hfill $\Box$ \\

Si les $f_{i}$ sont solutions d'une \'equation aux 
$q$-diff\'erences, nous dirons simplement dans ce cas 
que ces solutions sont ind\'ependantes. Notre but est
de construire une famille maximale de solutions
ind\'ependantes de l'\'equation (1). \\

\textsl{A.3.3 Lemme. -}
\emph{Le nombre de solutions ind\'ependantes de (1)
ne peut exc\'eder $n$, l'ordre de l'\'equation.} \\

\Pr
Soient en effet $f_{1},\ldots,f_{n+1}$ des solutions
de (1). Les lignes 
$L_{i} = (\sigma_{q}^{i} f_{1},\ldots,\sigma_{q}^{i}f_{n+1})$
sont alors li\'ees par la relation
$a_{0} L_{n} + \cdots + a_{n} L_{0}$ et l'on conclut
gr\^ace au lemme A.3.2.
\hfill $\Box$ \\

\textsl{A.3.4 Th\'eor\`eme. -}
\emph{Dans le cas formel, on peut construire $n$ solutions 
ind\'ependantes.} \\

\Pr
Elle se fait par r\'ecurrence sur l'ordre de 
l'op\'erateur $P$; l'algorithme correspondant est
r\'ecursif. On exploite naturellement les r\'esultats 
sur la factorisation de 1.2 et ceux sur les \'equations
du premier ordre avec second membre de A.1. \\

Si $n = 1$, on peut \'ecrire $P = a(z^{-\mu} \sigma - c) u^{-1}$,
et $u e_{q,c} \Theta_{q}^{\mu}$ est une solution non nulle. \\

Si $P$ est d'ordre $n = m + 1 \geq 2$, on \'ecrit
$P = a(z^{-\mu} \sigma - c) u^{-1} Q$, o\`u $Q$
est d'ordre $m$. Par hypoth\`ese de r\'ecurrence, 
il y a $m$ solutions ind\'ependantes $f_{1},\ldots,f_{m}$
de $Q$. D'apr\`es le th\'eor\`eme A.2.4, il existe
$f \in L$ tel que $Q f = u e_{q,c} \Theta_{q}^{\mu}$.
Il est clair que $f,f+f_{1},\ldots,f+f_{m}$ sont solutions
de $P$. \\

Par multilin\'earit\'e altern\'ee du d\'eterminant
le $q$-Wronskien de $f,f+f_{1},\ldots,f+f_{m}$ est
\'egal \`a celui de $f,f_{1},\ldots,f_{m}$. On
manipule les lignes de ce dernier on remplace $L_{m}$
par $b_{0} L_{m} + \cdots + b_{m} L_{0}$, 0\`u
$Q = b_{0} \sigma^{m} + \cdots + b_{m}$. Cela
multiplie le d\'eterminant par $b_{0}$. Mais cela
remplace aussi la derni\`ere ligne par 
$(Qf,Qf_{1},\ldots,Qf_{m}) = (Qf,0,\ldots,0)$.
Le coefficient $Qf$ vaut $u e_{q,c} \Theta_{q}^{\mu}$,
qui est inversible, et son cofacteur est le $q$-wronskien 
de $(f_{1},\ldots,f_{m})$. On obtient ainsi la formule:
$$
W_{q}(f,f+f_{1},\ldots,f+f_{m}) = 
\frac{1}{b_{0}} u e_{q,c} \Theta_{q}^{\mu} W_{q}(f,f_{1},\ldots,f_{m}).
$$
Il est donc non nul, ce qui ach\`eve la preuve.
\hfill $\Box$ 

\subsection*{A.4~~~R\'esolution analytique}

On se place ici dans le cas convergent. Si l'on reprend 
la factorisation $P = a(z^{-\mu} \sigma - c) u^{-1} Q$ 
exploit\'ee en A.3, on constate que l'on n'a la garantie
d'une factorisation convergente que si toutes les pentes
de $Q$ sont $\leq \mu$ (1.2.8). Mais, si l'une d'elles
est $\lt \mu$, le th\'eor\`eme A.2.4 ne s'applique pas.
Ainsi, la m\'ethode de A.3 ne s'applique \`a la 
r\'esolution convergente que si $S(P) = \{\mu\}$,
autrement dit, si $P$ est \emph{pur}.  On ne peut donc
esp\'erer trouver $n$ solutions ind\'ependantes en g\'en\'eral. \\

\textsl{A.4.1 Th\'eor\`eme (lemme d'Adams). -}
\emph{Soit $\mu_{k}$ la premi\`ere pente de $P$.
L'\'equation (1) admet alors $r_{P}(\mu_{k})$ 
solutions convergentes ind\'ependantes.} \\

\Pr
On d\'eduit en effet de A.3 une factorisation $P = QR$
avec $R$ pur de pente $\mu_{k}$ et d'ordre $r_{P}(\mu_{k})$.
On applique alors \`a $R$ la m\'ethode de A.3 (on est dans
la cas (i) du th\'eor\`eme A.2.4).
\hfill $\Box$

%%%%%%%%%%%%%%%%%%%%%%%%%%%%%%%%%%%%%%%%%%%%%%%%%%%%%%%%%%%%%%%%%%%%%%%%%%%%%

\end{document}